\documentclass[a4paper]{article}
\usepackage{geometry}
\usepackage{tabularx}
\usepackage{tikz}
\usepackage{subcaption}
\usepackage{rotating}
\usepackage{overpic}
\usepackage{amsmath, amssymb, amsthm,  amsfonts, mathrsfs, cite}
\usepackage{multirow}
\usepackage{color}
\usepackage{booktabs}
\usepackage{bm}
\usepackage{scrextend}
\usepackage{arydshln}
\usepackage {url}
\usepackage[font=scriptsize, labelfont=bf]{caption}
\usepackage{graphicx} 
\usepackage{epstopdf} 
\usepackage{algpseudocode}
\usepackage{array}
\usepackage{algorithm, algcompatible}
\usepackage{paralist} 
\usepackage{verbatim} 

\usepackage{xcolor, marginnote, enumitem}
\numberwithin{equation}{section}
\newtheorem{lemma}{Lemma}
\newtheorem{theorem}{Theorem}
\newtheorem{proposition}{Proposition}

\newtheorem{assumption}{Assumption}
\theoremstyle{remark}
\newtheorem{remark}{Remark}
\theoremstyle{definition}
\newtheorem{definition}{Definition}

\DeclareMathOperator{\dist}{dist}

\DeclareMathOperator{\dom}{dom}

\newcommand{\bE}{\mathbb{E}}
\newcommand{\bD}{\mathbb{D}}

\newcommand{\bitem}{\begin{itemize}}
\newcommand{\eitem}{\end{itemize}}

\newcommand{\bpm}{\begin{pmatrix}}

\DeclareMathAlphabet{\mathbfit}{OML}{cmm}{b}{it}

\usepackage[plainpages=false, colorlinks=true, 
            linkcolor=black, urlcolor=black, citecolor=black]{hyperref}

\begin{document}
\title{A boosted second-order convex splitting algorithm based on gradient flows}
\author{Xinhua Shen\thanks{School of Mathematics, 
		Renmin University of China,  China  \  \href{mailto:shenxinhua@ruc.edu.cn}{shenxinhua@ruc.edu.cn}}
\and Zaijiu Shang\thanks{1. Center for Mathematics and Interdisciplinary Sciences at Fudan University, Shanghai, 200433, China \  \href{mailto:zaijiu@amss.ac.cn}{zaijiu@amss.ac.cn} 2. Shanghai Institute for Mathematics and Interdisciplinary Sciences, Shanghai, 200433, China \ \href{mailto:zaijiu@simis.cn}{zaijiu@simis.cn}}
\and   Hongpeng Sun\thanks{School of Mathematics,
Renmin University of China, China \  \href{mailto:hpsun@amss.ac.cn}{hpsun@amss.ac.cn} }
}
\maketitle
\begin{abstract}
This paper introduces a second-order convex splitting scheme for gradient flows arising in phase-field models, based on the backward differentiation formula (BDF2) for the implicit part and the Adams-Bashforth method for the nonlinear and explicit component. The method is formulated and analyzed in finite-dimensional spaces, where energy stability plays a central role in establishing rigorous convergence properties. By leveraging the Kurdyka-\L ojasiewicz framework, we prove the global convergence of the discrete trajectories generated by the scheme, even in the presence of nonsmooth energy functionals, under mild assumptions on the time-step size. The Armijo line search strategy and the classical preconditioning strategies, such as symmetric Gauss-Seidel and Jacobi, are incorporated to improve its computational efficiency. Numerical experiments confirm that the proposed method achieves computational efficiency compared to existing first-order splitting approaches and other accelerated splitting algorithms, while maintaining robustness in both smooth and nonsmooth regimes.
\end{abstract}

\paragraph{MSCodes.}{65K10,
65F08,
49K35,
90C25.
}

\section{Introduction}\label{sec:intro}
We consider the following nonconvex optimization problem
\begin{equation}\label{eq:basefunctional}
    \min_{u \in X} E(u)= H(u) + F(u)
\end{equation}
where $X$ is a finite-dimensional Hilbert space, $H: X\to\mathbb{R}\cup\{+\infty\}$ is a proper, lower semicontinuous and convex function, $F$ is a proper function whose gradient $\nabla F$ is Lipschitz continuous with Lipschitz constant $L$. Nonconvex problems have wide-ranging applications, and the difference of convex functions algorithm (shortened as DCA henceforth) is one of the most efficient schemes \cite{LeThi2018, LPH, Lethi2024}. We refer to \cite{APX, Pangcui} for detailed discussions of the reader to the difference-of-convex-functions (DC) structure and analysis for many widely used nonconvex optimization problems.

Numerous works \cite{Schropp2000,helmke1996,su2016} have explored the connection between ordinary differential equations and optimization. In particular, the study \cite{su2016} offers novel insights into Nesterov's accelerated method and establishes a series of theoretical results regarding the first-order gradient method. Motivated by these investigations, the gradient descent method can be interpreted as discretizing the gradient flow $\frac{du}{dt}=-\nabla Q(u)$ via the forward Euler method, which takes the form
\begin{equation*}
    \frac{u^{n+1}-u^n}{\delta t} = -\nabla Q(u^n) ,\quad u^0=u(0)
\end{equation*}
where $Q$ is a smooth objective function and $\delta t>0$ is the time step.

We consider the subgradient dynamical system for \eqref{eq:basefunctional}
\begin{equation}
    \frac{du(t)}{dt} \in - \partial H(u(t)) - \nabla F(u(t)), \quad u(0) \in \dom E, \ t \in (0,\infty).
\end{equation}
While $u$ is an absolutely continuous curve, it is called a subgradient curve \cite{bolte1}.
The second-order implicit-explicit (IMEX) splitting technique is widely used in phase-field simulations \cite{Shen2010, Li2017, Feng2013} and for solving systems of ordinary differential equations \cite{ARW1995}. 
Inspired by the relationship between the Euler scheme and the gradient descent method, we develop a boosted algorithm for general nonsmooth and nonconvex problems. This algorithm discretizes the subgradient dynamic system (flow) using numerical discretization techniques based on the IMEX splitting method \cite{ARW1995}.
This yields the implicit scheme
\begin{equation}\label{eq:bdf2:adam:orig}
    \mathbf{0}\in\frac{2}{3\delta t}(3u^{n+1} - 4u^n + u^{n-1}) + \partial H(u^{n+1})+ 2f(u^n) - f(u^{n-1}),
\end{equation}
where $f=\nabla F$. The scheme reduces to the classical second-order IMEX splitting method, as discussed in \cite{ARW1995}, when $H$ is smooth with a Lipschitz-continuous gradient $\nabla H$.

It is well known that \eqref{eq:bdf2:adam:orig} is based on the second-order backward differentiation formulas (BDF2) for the implicit component $\partial H$ and the Adams-Bashforth scheme for the explicit component $f$ \cite{ARW1995, Shen2010}. Under different smoothness assumptions and line search strategies, several boosted DC methods \cite{Artacho2018, ArtachoSIAM, Ferreira2025} have been developed to improve the performance of the DC algorithm. The inexact and boosted  DC algorithm is developed in \cite{Ferreira2025} under nonsmooth assumptions via a nonmonotone line search strategy. Inspired by these boosted DC approaches with line search and recent advances in preconditioned frameworks for DC algorithms \cite{DS, Shensun2023}, we propose a preconditioned second-order convex splitting framework. It can be seen as a generalized DC with varying DC components, integrated with a line search scheme to solve  \eqref{eq:basefunctional}.

We introduce the following energy functions for the construction of our algorithm and the subsequent theoretical analysis:
\begin{equation}\label{eq:def:en:hn}
\begin{aligned}
    &{E}^n(u) = {H}^n(u) - {F}^n(u),\\
    &{H}^n(u) = H(u) +\frac{1}{\delta t} \|u-u^n\|^2,  \\
    &{F}^n(u) = \frac{1}{3\delta t}\|u-u^{n-1}\|^2 - F(u)-\langle f(u^n) - f(u^{n-1}),u - u^{n-1}\rangle.
\end{aligned}
\end{equation}
We subsequently solve the subproblem involving the difference of convex functions by incorporating a proximal term to compute the update $y^n$.
 The optimization problem for determining $y^n$ is formulated as follows:
\begin{equation}\label{eq:mini:dc:sub}
    y^n = \arg\min_u \left\{{H}^n(u) - \langle \nabla {F}^n(u^n), u\rangle +\frac12\|u-\hat{u}^n\|_M^2\right\}.
\end{equation}
Here, the positive-semidefinite and self-adjoint operator (metric) $M$ is used to construct efficient preconditioners, with $\|u\|_{M}^2:= \langle Mu, u\rangle$. The starting point for preconditioned iteration is denoted by $\hat{u}^n$, which can be chosen as either $\hat{u}^n= u^n$ or $\hat{u}^n= \tilde u^n=\frac{4}{3}u^n - \frac{1}{3}u^{n-1}$. 
The initial guess $\hat u^n$ admits two choices in two ways.  The first choice, $\hat u^n = u^n$, corresponds to the standard proximal method. 
The second choice is $\hat u^n = \tilde u^n$, where
\begin{equation*}
    \tilde u^n = \frac{4}{3}u^n - \frac{1}{3}u^{n-1} = u^n + \frac{1}{3}\bigl(u^n - u^{n-1}\bigr).
\end{equation*}
It can be interpreted as an extrapolation step with a fixed extrapolation parameter $1/3$ at $u^n$. A further motivation for introducing $\tilde u^n$ stems from the discretization scheme.  It can be verified that $u^{n+1} - \tilde u^n = \frac{1}{3}\bigl(3u^{n+1} - 4u^n + u^{n-1}\bigr)$, which is precisely the discretization of the second-order backward differentiation formula (BDF2).  

Applying the first-order optimality condition, the minimization problem \eqref{eq:mini:dc:sub} for determining $y^n$ yields the equation:
\begin{equation}\label{eq:solve:y:pre:general}
    \frac{2}{3\delta t}(3y^n - 4u^n + u^{n-1}) + M(y^n - \hat{u}^n) + v^n + (2f(u^n) - f(u^{n-1})) = \mathbf{0},
\end{equation}
where $v^n\in \partial H(y^n)$. We reformulate \eqref{eq:solve:y:pre:general} as classical preconditioned iterations for solving $y^n$, which simplifies the solution of $y^n$, offering greater convenience. For further details on preconditioning techniques, see Section \ref{sec:pre:sub}. It is crucial to highlight that the inclusion of the extrapolated gradient term $2f(u^n) - f(u^{n-1})$, distinguishes the scheme \eqref{eq:bdf2:adam:orig} substantially from DCA-like algorithms \cite{Lethi2021}. 
For the Allen-Cahn model with specific $H$ and $F$ functions, where $S$ represents a positive constant and $I$ is the identity operator, the choices discussed in \cite{Shen2010} involve $M=S I$ and $\hat u^n= 2u^n-u^{n-1},$ while the corresponding settings in \cite{Li2017} consider $M = \delta t S I$ and $\hat u^n = u^n$. The careful selection of $S(u^{n+1}- 2u^n + u^{n-1})$ in \cite{Shen2010} or $S\delta t (u^{n+1}-u^n)$ in \cite{Li2017} preserves second-order time discretization and improves numerical stability. Incorporating additional information into the iteration process is a widely adopted technique in optimization and numerical methods. The reflected gradient method \cite{Malitsky2015, Tseng2000} utilizes the value at the reflection point, expressed as $2u^n - u^{n-1}$, to yield enhanced interpretability of the objective function. We incorporate not only the prior information from the iterative sequence $\{u^n\}_n$, such as $3y^{n+1} - 4u^n + u^{n-1}$, but also the corresponding gradient component, represented by $2f(u^n) - f(u^{n-1})$.

By defining $d^n = y^n-u^n$, we proceed with an $\emph{aggressive}$ Armijo-type line search involving parameters $\beta^k \lambda_0>0$, $k = 0,1,2, \ldots$, $\beta \in (0,1)$, and $\lambda_0, \alpha>0$. The objective is to find the smallest integer $k$ satisfying:
\begin{equation}\label{eq:line:search}
\hat E(y^n+\beta^k \lambda_0 d^n) \leq \hat E(y^n) - \alpha\beta^k \lambda_0 \|d^n\|^2, \quad d^n  = y^n-u^n. 
\end{equation}
Here $\hat E$ represents the line search energy. The final step size $\lambda_n := \beta^k \lambda_0$ obtained from the line search in \eqref{eq:line:search} dictates the update of $u^{n+1}$ as follows:
\begin{equation}\label{eq:LSstep}
   u^{n+1} = y^n + {\lambda}_n d^n.
\end{equation}
In cases where the line search is viable, the relationship $u^{n+1} = u^n + (1+\lambda_n)d^n$ holds, with $\lambda_n >0$. This signifies that the update $u^{n+1}$ involves a full step of $d^n$ alongside an additional step of $\lambda_n d^n$. This approach differs from the Fukushima-Mine line search \cite{Mine1981}, as well as traditional line search methods used in Newton techniques \cite[Algorithm 2.1]{Qi1999}, and the Armijo line search employed in generalized conditional gradient methods \cite{Karl}, where typically $(1+\lambda_n) \in  [0,1]$ with $\lambda_n\in[-1,0]$.

Our research contributions can be summarized as follows. First, we establish the global convergence of the iterative sequence generated by the second-order convex splitting method. Second-order convex splitting methods have been widely applied to phase-field simulations involving the Allen–Cahn, Cahn–Hilliard equations, and gradient flow problems \cite{ARW1995, Feng2013, Li2017, Shen2010}. Existing convergence analyses have focused primarily on energy stability or perturbed energy stability \cite[Theorem 3.2]{Feng2013}, \cite[Lemma 2.3]{Shen2010}, or \cite[Theorem 1.2]{Li2017}. Recent advancements in Kurdyka-\L ojasiewicz (KL) analysis  \cite{Attouch2009, ABS, Li2018} allow us to prove global convergence of the iterative sequence under mild assumptions. In addition, we integrate line search acceleration \cite{Artacho2018, ArtachoSIAM, Artacho2022} to accelerate the second-order convex splitting method and incorporate a preconditioning technique \cite{DS, Shensun2023} to solve large-scale linear subproblems efficiently at each iteration. By performing a finite number of preconditioned iterations without error control, we ensure global convergence while retaining the benefits of line-search acceleration. Numerical experiments verify the superior computational efficiency of the proposed preconditioned second-order convex splitting method integrated with line search.

The remainder of this paper is organized as follows. Section \ref{sec:conv:pre_ls_1} introduces some essential analytical tools and presents the proposed second-order convex splitting scheme, which incorporates line search and preconditioning techniques. A comprehensive global convergence analysis with KL properties under different assumptions is introduced in Section \ref{sec:conv_analysis}. In Section \ref{sec:pre}, we discuss the preconditioning technique and present detailed numerical tests to demonstrate the effectiveness of the proposed splitting scheme. Finally, we discuss the implications of the preconditioning and line search strategies and summarize our main findings in Section \ref{sec:conclusion}.
\section{Algorithm and convergence analysis}\label{sec:conv:pre_ls_1}
\subsection{Preliminary}
For a proper, lower semicontinuous (closed) function $ g : \mathbb{R}^n \to \mathbb{R} \cup \{\infty\} $, its (limiting) subdifferential of $ g $ at $ x \in \text{dom } g $ is defined as
\begin{align}\label{def:lim_subdifferential}
    \partial g(x) = \left\{ v \in \mathbb{R}^n : \exists\ x^t \xrightarrow{g} x, v^t \to v,\liminf_{y \to x^t} \frac{g(y) - g(x^t) - \langle v^t, y - x^t \rangle}{\|y - x^t\|} \geq 0,\ \forall t\right\},
\end{align}
where $ z \xrightarrow{g} x $ denotes that $ z \to x $ and $ g(z) \to g(x) $. The above subdifferential reduces to the classical subdifferential in convex analysis when $ g $ is convex, i.e.,
\begin{equation*}
    \partial g(x) = \left\{ v \in \mathbb{R}^n : g(u) - g(x) - \langle v, u - x \rangle \geq 0, \forall u \in \mathbb{R}^n \right\}.
\end{equation*}
In addition, if $ g $ is continuously differentiable, then the subdifferential \eqref{def:lim_subdifferential} reduces to the gradient of $ g $, denoted by $ \nabla g $. For additional details, please refer to \cite{Roc1, Wen2018}.

We introduce the following Kurdyka-\L ojasiewicz (KL) property and KL exponent for the global and local convergence analysis. While the KL properties \cite[Definition 2.4]{ABS} \cite[Definition 1]{Artacho2018} enable establishing the global convergence of iterative sequences, the KL exponent \cite[Remark 6]{Bolte2014} assists in determining a local convergence rate. 
\begin{definition}[KL property, KL function and KL exponent]\label{def:KL}Let $g: \mathbb{R}^n \rightarrow \mathbb{R}$ be a proper, closed function. $g$ is said to satisfy the KL property if for any critical point $\bar x$,  there exists $\nu \in (0,+\infty]$, a neighborhood $\mathcal{O}$ of $\bar x$, and a continuous concave function $\psi: [0,\nu) \rightarrow [0,+\infty)$ with $\psi(0)=0$ such that:
	\begin{itemize}
		\item [{\emph{(i)}}] $\psi$ is continuously differentiable on $(0,\nu)$ with $\psi'>0$ over $(0,\nu)$;
		\item [{\emph{(ii)}}] for any $x \in \mathcal{O}$ with $g(\bar x) < g(x) <g(\bar x) + \nu$, one has
		\begin{equation}\label{eq:kl:def}
		\psi'(g(x)-g(\bar x)) \cdot {\dist(\mathbf{0},\partial g(x))}\geq 1. 
		\end{equation}
	\end{itemize}
	Furthermore, for a proper, closed function $g$ satisfying the KL property,  if $\psi$ in \eqref{eq:kl:def} can be chosen as $\psi(s) = cs^{1-\theta}$ for some $\theta \in [0,1)$ and $c>0$, i.e., there exist  $\bar c, \epsilon >0$ such that
	\begin{equation}\label{eq:KL:exponent:theta:exam}
	{\dist(\mathbf{{0}},\partial g(x))}\geq \bar c (g(x)-g(\bar x))^{\theta}
	\end{equation}
 whenever $\|x -\bar x\| \leq \epsilon$ and $g(\bar x) < g(x) <g(\bar x) + \nu$, then we say that $g$ has the KL property at $\bar x$ with exponent $\theta$.
\end{definition}
The KL exponent is determined solely by the critical points. If $g$ exhibits the KL property with an exponent $\theta$ at any critical point $\bar x$, it follows that $g$ is a KL function with an exponent $\theta$ at all points within the domain of definition of $\partial  g$ \cite[Lemma 2.1]{Li2018}.

A useful conclusion is that a semi-algebraic function (see \cite[Definition 2.1]{ABS} for the definitions of semi-algebraic sets and semi-algebraic functions) has the KL property and is a KL function \cite[Lemma 2.3]{ABS}. Additionally, we assume that the energy $E(x)$ is level-bounded, defined by $\text{lev}_{\leq \alpha}(E): = \{x: E(x) \leq \alpha\}$ being bounded (or possibly empty) \cite[Definition 1.8]{Roc1}. It is worth noting that a function $f: \mathbb{R}^n \rightarrow \mathbb{R}$ being coercive (i.e., $f(x) \rightarrow +\infty$ as $\|x\|\rightarrow +\infty$) implies that it is also level-bounded.  The following lemma facilitates the analysis of global convergence.
\begin{lemma}\label{lem:level_bound_A}
    Let $x, y \in X$, and suppose $E(x)$ is level-bounded. Then
\begin{equation*}
    A(x,y) = E(x) + \gamma \|x - y\|_{\mathcal{M}}^2
\end{equation*}
is also level-bounded, where $\gamma > 0$ is a constant and $\mathcal{M}$ is positive definite. Moreover, if $E(x)$ is a Kurdyka--\L ojasiewicz (KL) function, then $A(x,y)$ is also a KL function.
\end{lemma}
\begin{proof}
    Let $ S_t := \{ (x,y) \in X\times X \mid A(x,y) \le t \} $. Since $E(x)$ is level-bounded, the boundedness of $ x $ can be established by noting that $ E(x) \le E(x) + \gamma \|x - y\|_{\mathcal{M}}^2\le t $. In addition, note that $E(x)$ is also proper, one obtains that there exists a constant $U_1>0$ that satisfies $\lambda_{\mathcal{M}}\|x-y\|^2\le\|x-y\|^2_\mathcal{M} \le \lambda_{\mathcal{M}}U_1^2$, where $\lambda_{\mathcal{M}}>0$ is the minimum eigenvalue of $\mathcal{M}$. With the boundedness of $x$, i.e., $\|x\|\le U_2, \forall (x,y)\in S_t$, the estimate $\|y\|\le\|y-x\| + \|x\|\le U_1 + U_2$ holds for any $(x,y) \in S_t$. By the definition of the product norm in Hilbert space, the proof of the level-boundedness of $A(x,y)$ is completed.

   To establish the KL property of $ A(x,y) $, we can conclude this by applying the existing theorem \cite[Lemma 4]{DS}, thereby completing the proof.
\end{proof}

\subsection{Algorithm and energy decay}
Before presenting our algorithm, we discuss the step size $\delta t$ introduced by the BDF scheme. The incorporation of the parameter $\delta t$ improves the flexibility of our algorithm. To ensure the convergence of our algorithm, it is necessary to impose a restriction on the parameter $\delta t$. We mainly need the convexity of $F^n$, and Lemma \ref{lem:con:L} guarantees this property via a suitable restriction on the parameter $\delta t$. 
\begin{lemma}\label{lem:con:L}
    Assume that the Lipschitz constant of $f$ is $L > 0$. When $\delta t < \frac{2}{3L}$, $F^n$ is strongly convex.
\end{lemma}
    \begin{proof}
        Direct calculation gives \[
        \nabla F^n(u) = \frac{2}{3\delta t}(u-u^{n-1}) - f(u) -(f(u^n)-f(u^{n-1})).
        \]
        A direct computation shows that for any $u_1, u_2 \in X$, the following holds
        \begin{align*}
           & \langle \nabla F^n(u_1)  - \nabla F^n(u_2), u_1-u_2\rangle = \langle  \frac{2}{3\delta t}(u_1 - u_2) -f(u_1) + f(u_2), u_1 - u_2\rangle  \\
           &= \frac{2}{3\delta t} \|u_1 - u_2\|_{2}^2 - \langle f(u_1) - f(u_2), u_1 - u_2  \rangle \geq ( \frac{2}{3\delta t} -L)\|u_1 - u_2\|_{2}^2. 
        \end{align*}
        We conclude that $\nabla F^n$ is strongly monotone with constant $( \frac{2}{3\delta t} -L)$ under the condition $\delta t < \frac{2}{3L}$, implying that $F^n$ is strongly convex. This completes the proof.
    \end{proof}
With the constraint on $\delta t$, we now present an outline of our algorithm.
\begin{algorithm}[H]
\renewcommand{\algorithmicrequire}{\textbf{Input:}}
 \renewcommand{\algorithmicensure}{\textbf{Output:}}
\caption{Based on gradient \textbf{flow}, a \textbf{b}oosted algorithm with \textbf{B}DF and \textbf{A}dams-Bashforth convex splitting using \textbf{p}reconditioning for subproblem (abbreviated as \protect{\textbf{FlowBBAp}})}\label{alg:cap} 
\begin{algorithmic}[1]
\State Choose $u^0$, $0<\delta t<2/(3L)$, $\alpha>0$, $0<\beta<1$, $\bar\lambda_{\max}>0$. Set $u^{-1}=u^0$
\State $H^n(u) := H(u) + \frac{1}{\delta t}\|u-u^n\|^2$
\State $F^n(u) := \frac{1}{3\delta t}\|u - u^{n-1}\|^2 - F(u) -\langle f(u^n) - f(u^{n-1}),u-u^{n-1}\rangle$
\State Solve the following subproblem 
\begin{equation}\label{eq:algorithmic:update}
    y^{n} = \arg\min_u \left\{H^n(u) - \langle \nabla F^n(u^{n}),u\rangle + \frac12\|u-\hat{u}\|^2_M\right\}
\end{equation}
where $\hat{u}$ can be chosen from the set $\{u^n, \tilde{u}^n\}$ throughout the algorithm, where $\tilde u^n=\frac{4}{3}u^{n}-\frac{1}{3}u^{n-1}$.
\State Set $d^n := y^n - u^n$. IF $d^n = 0$, STOP and RETURN $u^n$. Otherwise, go to Step 6.
\State line search process with $\alpha$, $d^n$, $\bar \lambda_{\max}$ and $\beta$: IF `Line Search Criterion' (\eqref{ineq:MS_LS_rule} for `MPCLS' or \eqref{ineq:ori_LS_rule} for NLS) is satisfied with line search step $\lambda_n$, set $u^{n+1}:= y_n + \lambda_nd^n$, set $n = n+1$ and go to Step 4.
\end{algorithmic}
\end{algorithm}

In Step 4 of Algorithm \ref{alg:cap}, there are two options for defining $\hat u$, each corresponding to a different algorithm: (1) if we set $\hat u = u^n$, this version is denoted FlowBBAp; (2) if we define $\hat u = \tilde u^n$, it is denoted FlowBBApe. 

Step 6 of Algorithm \ref{alg:cap} involves two different line search strategies: the `Multistep Predictive Compatible Line Search (MPCLS)' and the `Normal Line Search (NLS)'. Both methodologies rely on the potential search direction $d^n$ established in Step 5. A key advantage of the `MPCLS' strategy is that the direction $d^n$ serves as a valid descent direction for $E^n$ as established in Proposition \ref{prop:descent_2} later. The specifics of both strategies are detailed below.
\begin{itemize}
    \item \textbf{Multistep Predictive Compatible Line Search (MPCLS)}: Its `Line Search Criterion' is defined as 
    \begin{equation}\label{ineq:MS_LS_rule}
    E^n(y_n + \lambda_nd^n)\le E^n(y^n) - \alpha\lambda_n\|d^n\|^2.
    \end{equation}
    The procedure for Step 6 of Algorithm \ref{alg:cap} can be reformulated as follows:
    \begin{align*}
        &\text{Set } \lambda_n = \bar\lambda_{\max}.\ \text{WHILE } E^n(y_n + \lambda_nd^n)> E^n(y^n) - \alpha\lambda_n\|d^n\|^2, \text{ DO } \lambda_n = \beta\lambda_n.\\
        &\text{Set }u^{n+1} = y^n + \lambda_nd^n,\ n = n+1\ \text{and go to Step 4.}  
    \end{align*}

    \item \textbf{Normal Line Search (NLS)}: Its `Line Search Criterion' is defined as \begin{equation}\label{ineq:ori_LS_rule}
    E(y_n + \lambda_nd^n)\le E(y^n) - \alpha\lambda_n\|d^n\|^2.
    \end{equation}
    The procedure for Step 6 of Algorithm \ref{alg:cap} can be adjusted as follows:
    \begin{align*}
        &\text{IF}\ \langle \nabla E(y^n),d^n\rangle>0,\ \text{set $u^{n+1} =y^n,\ n = n+1$ and go to Step 4.}\\
        &\text{Set } \lambda_n = \bar\lambda_{\max}.\ \text{WHILE } E(y_n + \lambda_nd^n)> E(y^n) - \alpha\lambda_n\|d^n\|^2, \text{ DO } \lambda_n = \beta\lambda_n.\\
        &\text{Set }u^{n+1} = y^n + \lambda_nd^n,\ n = n+1\ \text{and go to Step 4.}  
    \end{align*}
\end{itemize}

For the line search step, we observe that the quadratic interpolation for finding an initial value of the line search step is highly beneficial; see the following remark.
\begin{remark}
    In Algorithm \ref{alg:cap}, instead of starting the line search process at point $u^n$ as done in the classical Armijo line search, we initiate the process at the newly computed point $y^n$ derived from \eqref{eq:algorithmic:update}. This modification aligns with Proposition \ref{lem:LSbase}, indicating that the transition from $y^n$ to the subsequent point $u^{n+1}$ exhibits an energy decay property. Selecting a more precise point close to the optimal solution along the line search direction, rather than a maximum searching step, i.e., $\lambda_n:= \bar \lambda_{\max}$, can save computational time. At the point $y^n$, a quadratic interpolation is employed to estimate the function $e_n(\lambda):= \hat E(y^n + \lambda d^n)$, where $\hat E$ selects either $E$ or $E^n$ based on varying line search criteria. By considering three key pieces of information: $e_n(0) = \hat E(y^k)$, $e_n'(0) = \langle \nabla \hat E(y^n), d^n\rangle $, and $e_n(\bar{\lambda}) =\hat E(y^n + \bar{\lambda} d^n)$, an approximation function can be formulated as \cite[Corollary 2]{Artacho2018} or \cite[Chapter 3.4]{NoceWrig06}
\begin{equation}
    \mathcal{E}_n(\lambda) = \frac{e_n(\bar{\lambda}) - e_n(0) - \bar{\lambda}e_n'(0)}{\bar{\lambda}^2}\lambda^2 + e_n'(0)\lambda + e_n(0).
\end{equation}
Moreover, the minimum point of the quadratic function $\mathcal{E}_n(\lambda)$ is located at
\begin{equation*}
    \lambda^*_n = \frac{e_n'(0)\bar{\lambda}^2}{2(e_n(\bar{\lambda}) - e_n(0) - e_n'(0)\bar{\lambda})}.
\end{equation*}
Hence, the algorithm's computational efficiency can be accelerated by replacing the setting $\lambda_n:=\bar \lambda_{\max}$ in Step 6 with $\lambda_n:=\min\{\bar \lambda_{\max},\lambda^*_n\}$.
\end{remark}

Our discussion, based on Algorithm \ref{alg:cap}, starts with the property of energy decay. The subsequent proposition characterizes the energy reduction from $u^n$ to $y^n$.
\begin{proposition}\label{prop:descent}
Let $\{y^n\}_n$ and $\{u^n\}_n$ be sequences generated iteratively by solving \eqref{eq:mini:dc:sub} and \eqref{eq:LSstep} with $\hat u^n=u^n$. Define $N_1 = (\frac{4}{3\delta t}  - \frac{L}{2})I + M$. Then it holds that $E^n(y^n) \leq  E^n(u^n) - \|d^n\|^2_{N_1}$.
\end{proposition}
\begin{proof}
    By the strong convexity of $H^n(u)$ and $F^n(u)$, the following inequalities hold
    \begin{align*}
        H^n(u^n) - H^n(y^n) &\geq \langle v^n,u^n - y^n\rangle + \frac{1}{\delta t}\|u^n - y^n\|^2, \\
        F^n(y^n) - F^n(u^n) &\geq \langle\nabla F^n(u^n),y^n - u^n\rangle + (\frac{1}{3\delta t} - \frac{L}{2})\|y^n - u^n\|^2,
    \end{align*}
    for any $v^n\in\partial H^n(y^n)$. Combining these two inequalities yields
    \begin{equation*}
    E^n(u^n) - E^n(y^n) \geq \langle v^n-\nabla F^n(u^n) ,u^n - y^n\rangle+ (\frac{4}{3\delta t} - \frac{L}{2})\|y^n - u^n\|^2.
    \end{equation*}
    From \eqref{eq:mini:dc:sub}, $y^n$ is a critical point of this convex function and satisfies the following optimality condition
    \begin{equation*}
    - M(y^n - u^n)\in\partial H^n(y^n) -\nabla F^n(u^n).
    \end{equation*}
    Combining these results, we arrive at
    \begin{equation*}
    E^n(u^n) - E^n(y^n) \geq (\frac{4}{3\delta t} - \frac{L}{2})\|y^n - u^n\|^2 + \|y^n-u^n\|_M^2.
    \end{equation*}
    This completes the proof.
\end{proof}

During an iteration, we have already observed a partial decrease in energy from $u^n$ to $y^n$, while the remaining energy decay can be derived from the line search process. The line search is not only essential for accelerations, but also for energy decay. In particular, we will show that $d^n$ serves as a descent direction of $E^n(u)$ on $y^n$ for the flexibility of the line search in the next section. By following a procedure analogous to the demonstration of Proposition \ref{prop:descent}, we can establish the connection between $E^n(y^n)$ and $E^n(u^n)$, where $\hat u^n = \tilde u^n$.

\begin{proposition}\label{prop:til:1}
Let $\{y^n\}_n$ and $\{u^n\}_n$ be sequences generated iteratively by solving \eqref{eq:mini:dc:sub} and \eqref{eq:LSstep} with $\hat u^n=\tilde u^n$. Then the following estimate holds:
    \begin{equation}\label{ineq:pre:prop1}
    E^n(y^n) \leq  E^n(u^n) - ( \frac{4}{3\delta t} - \frac{L}{2})\|y^n - u^n\|^2 - \frac{5}{6}\|y^n - u^n\|^2_M + \frac{1}{6}\|u^n-u^{n-1}\|^2_M.
\end{equation}
\end{proposition}

\section{The global convergence for the second-order convex splitting scheme}\label{sec:conv_analysis}
In this section, we will investigate the global convergence of the proposed algorithmic frameworks. We first address the algorithm FlowBBAp with the `MPCLS' line search criterion, which is analyzed under a mild smoothness assumption. Next, we analyze the algorithm FlowBBApe with the `NLS' line search criterion, based on the assumption introduced in Section \ref{sec:intro}.
\subsection{Convergence analysis of FlowBBAp in the smooth setting}\label{sec:global_converg_un}
In this subsection, the global convergence analysis relies on the following additional assumption.
\begin{assumption}\label{assum:1} 
Assume that $H$ is a smooth function, where $h=\nabla H$ is a Lipschitz continuous function with Lipschitz constant $L_h$.
\end{assumption}

Under the smoothness assumption, the line search direction $d^n$, derived from the corresponding subproblem \eqref{eq:algorithmic:update}, can be theoretically guaranteed to remain a direction of energy descent for our proposed energy $E^n$. Moreover, a significant energy reduction can be achieved in this direction. Proposition \ref{prop:descent_2} states that the direction $d^n$ remains a descent direction. Lemma \ref{lem:LSbase} provides an estimate of the energy decay in the line search process.
\begin{proposition}\label{prop:descent_2}
Let $\{y^n\}_n$ and $\{u^n\}_n$ be sequences generated iteratively by solving \eqref{eq:mini:dc:sub} and \eqref{eq:LSstep} with $\hat u^n=u^n$. Then it holds that $\langle\nabla E^n(y^n), d^n\rangle\leq - \|d^n\|^2_{N_2}$, where $N_2 = (\frac{2}{3\delta t} - L)I + M$.
\end{proposition}
\begin{proof}
    We first evaluate the gradient of $E^n(u)$ at $y^n$
    \begin{equation}\label{eq:gradEn}
        \nabla E^n(y^n) = h(y^n) + \frac{2}{\delta t}(y^n - u^n) - \frac{2}{3\delta t}(y^n - u^{n-1}) + f(y^n) + f(u^n) - f(u^{n-1}).
    \end{equation}
    Combining this with the first optimality condition,
    \begin{equation}\label{eq:subpro:optcon}
    \nabla H^n(y^n) -\nabla F^n(u^n) + M(y^n - u^n) = \mathbf{0}.
    \end{equation}
    We rewrite the gradient $ \nabla E^n(y^n)$ into
    \begin{equation*}
    \nabla E^n(y^n) = \frac{2}{3\delta t}(u^n - y^n) - f(u^n) + f(y^n) - M(y^n - u^n).
    \end{equation*}
    We thus have
    \begin{align*}
    \langle\nabla E^n(y^n),d^n\rangle = &\langle\frac{2}{3\delta t}(u^n - y^n) - f(u^n) + f(y^n) - M(y^n - u^n), y^n - u^n\rangle\\
    =&-\frac{2}{3\delta t}\|y^n - u^n\|^2- \|y^n - u^n\|_M^2 + \langle f(y^n) - f(u^n),y^n-u^n\rangle\\
    \leq&-\frac{2}{3\delta t}\|y^n - u^n\|^2 + L\|y^n - u^n\|^2 -\|y^n - u^n\|_M^2  \\
    =& (L - \frac{2}{3\delta t} )\|y^n - u^n\|^2- \|y^n - u^n\|_M^2 
    \end{align*}
    which leads to the conclusion. 
\end{proof}

With the Proposition \ref{prop:descent_2}, we can establish the feasibility of the line search with the following lemma. 
\begin{lemma}\label{lem:LSbase}
Let $\{(y^n, \lambda_n, d^n)\}_n$ be generated by Algorithm \ref{alg:cap} with $\hat u^n=u^n$ and line search criterion \eqref{ineq:MS_LS_rule}. There exist a small positive parameter $\alpha$ and a positive line search step $\lambda_n $ such that 
\begin{equation}\label{eq:line:search:2}
    E^n(y^n + \lambda_n d^n) - E^n(y^n)  \leq -\alpha\lambda_n\|d^n\|^2.
\end{equation}
\end{lemma}
\begin{proof}
If $d^n = 0$, there is nothing to prove. Otherwise, applying the mean value theorem to $ E^n(y^n + \lambda d^n)$ regarding $\lambda$, there is some $t_{\lambda}\in (0,1)$ such that
\begin{align*}
    E^n(y^n + \lambda d^n) - E^n(y^n) &= \langle\nabla E^n(y^n+t_{\lambda}\lambda d^n),\lambda d^n\rangle\\
    & = \lambda\langle\nabla E^n(y^n),d^n\rangle + \lambda\langle\nabla E^n(y^n + t_{\lambda}\lambda d^n) - \nabla E^n(y^n),d^n\rangle\\
    & \leq -\lambda\|d^n\|^2_{N_2} + \lambda \| \nabla E^n(y^n + t_{\lambda}\lambda d^n) - \nabla E^n(y^n)\| \|d^n\|.
\end{align*}
As $\nabla E^n$ is continuous at $y^n$, there are some small $\delta, \alpha > 0$ with $\beta = \frac{2}{3\delta t} - L$ and $\alpha < \beta$, such that
\begin{equation*}
    \|\nabla E^n(z) - \nabla E^n(y^n)\|\leq(\beta - \alpha)\|d^n\| \ \text{whatever} \ \|z-y^n\|\leq\delta.
\end{equation*}
Since $\|y^n+t_{\lambda}\lambda d^n-y^n\|=t_{\lambda}\lambda\|d^n\|\leq\lambda \|d^n\|$, then for all $\lambda\in(0,\frac{\delta}{\|d^n\|})$, we deduce 
\begin{equation*}
    E^n(y^n + \lambda d^n) - E^n(y^n) \leq -\lambda\|d^n\|^2_{N_2} + (\beta - \alpha)\lambda\|d^n\|^2 \leq -\alpha\lambda\|d^n\|^2,
\end{equation*}
with Proposition \ref{prop:descent_2} and this proof is complete by taking $\lambda_n=\lambda$. 
\end{proof}

Lemma \ref{lem:LSbase} and Proposition \ref{prop:descent} describe the complete energy reduction mechanism in a single iteration. This decrease in energy leads to an initial convergence analysis, as shown in the following lemma. Furthermore, this analysis acts as a preliminary step towards establishing the global convergence of the sequence $\{u^n\}_n$. We denote $\lambda_{\min} : =\min\{ \lambda_n, \ n \geq 1\}$ and $\lambda_{\max} : =\max\{ \lambda_n, \ n \geq 1\}$ throughout this paper. 
\begin{lemma}\label{lem:square:summable}
Let $\{u^n\}_n$ be a sequence generated by Algorithm \ref{alg:cap} with $\hat u^n=u^n$ and line search criterion \eqref{ineq:MS_LS_rule} for solving \eqref{eq:basefunctional}. If  the line search step satisfies $\lambda_{\max}<\sqrt{\frac{4}{3\delta tL} - \frac12}-1$, the sequence $\{ \|u^{n+1}-u^n\|\}_n$ is square summable, i.e., 
    \begin{equation}\label{ineq:quadconv}
    \sum_{n=1}^{\infty}\|u^{n+1} - u^n\|^2 < \infty.
\end{equation}
\end{lemma}
\begin{proof}
First, combining Proposition \ref{prop:descent} with the Lemma \ref{lem:LSbase} and remembering $u^{n+1} = y^n + \lambda_n d^n$, we can derive a lower bound of $E^n(u^n) - E^n(u^{n+1})$ 
\begin{equation}\label{eq:lower:bound:en}
    \alpha\lambda_n\|d^n\|^2 +(\frac{4}{3\delta t} - \frac{L}{2}) \|d^n\|^2\leq\alpha\lambda_n\|d^n\|^2 + \|d^n\|^2_{N_1} \leq E^n(u^n) - E^n(u^{n+1}).
\end{equation}
Second, let us turn to the upper bound of $E^n(u^n) - E^n(u^{n+1})$. From the definition of $E^n(u)$, we can compute $E^n(u^n) - E^n(u^{n+1})$ as follows
\begin{align}
    &E^n(u^n) - E^n(u^{n+1}) \notag \\
    =&H(u^n) - H(u^{n+1}) + F(u^n) - F(u^{n+1}) + \langle f(u^n) - f(u^{n-1}),u^n - u^{n+1}\rangle \notag \\
    -& \frac{1}{3\delta t}\|u^n-u^{n-1}\|^2-\frac{1}{\delta t} \|u^{n+1}-u^n\|^2 + \frac{1}{3\delta t}\|u^{n+1}-u^{n-1}\|^2 \notag  \\
    \leq &E(u^n)- E(u^{n+1}) + \frac{L}{2}(\|u^n - u^{n-1}\|^2 + \|u^n - u^{n+1}\|^2) - \frac{1}{3\delta t}\|u^n-u^{n-1}\|^2 \notag  \\ 
    -& \frac{1}{\delta t} \|u^{n+1}-u^n\|^2 + \frac{2}{3\delta t}(\|u^n - u^{n-1}\|^2 + \|u^n - u^{n+1}\|^2) \notag  \\
    \leq&E(u^n)- E(u^{n+1}) + (\frac{L}{2} + \frac{1}{3\delta t})\|u^n - u^{n-1}\|^2 - (\frac{L}{2} + \frac{1}{3\delta t})\|u^{n+1} - u^n\|^2 \notag \\
    + &L\|u^{n+1} - u^n\|^2, \label{eq:esti:en:en+1}
\end{align}
where the first inequality uses the following Lipschitz continuity of $f(u)$, and the inequality $\|u^{n+1}-u^{n-1}\|^2 \leq 2(\|u^{n+1}-u^{n}\|^2 + \|u^{n}-u^{n-1}\|^2)$,
$$ \langle f(u^n) - f(u^{n-1}),u^n - u^{n+1}\rangle\leq \frac{L}{2}(\|u^n - u^{n-1}\|^2 + \|u^n - u^{n+1}\|^2).$$
Now, we define an auxiliary function $A(x,y) = E(x) + (\frac{L}{2}+\frac{1}{3\delta t})\|x-y\|^2$, then the estimate \eqref{eq:esti:en:en+1} can be rearranged as
\begin{equation}
    E^n(u^n) - E^n(u^{n+1}) \leq A(u^n,u^{n-1}) - A(u^{n+1},u^n) + L\|u^{n+1}-u^n\|^2.
\end{equation}
Furthermore, making use of the fact $d^n = y^n-u^n = \frac{1}{1+\lambda_n}(u^{n+1} - u^n)$ and the lower bound estimate \eqref{eq:lower:bound:en}, we finally arrive at
\begin{equation}\label{eq:esti:vary:lambdan}
    \left(\big(\alpha\lambda_n+ \frac{4}{3\delta t} - \frac{L}{2}\big) \frac{1}{(1+\lambda_n)^2}  - L\right)\|u^{n+1} - u^n\|^2\leq A(u^n,u^{n-1}) - A(u^{n+1},u^n).
\end{equation}
Since $0\le\lambda_n< \sqrt{\frac{4}{3\delta tL} - \frac12}-1$, together with $1/L > 3\delta t/2$ by Lemma \ref{lem:con:L}, it can be readily verified that the coefficient $c(\lambda_n, \delta t, L) := \frac{4}{3\delta t} - \frac{L}{2} - L(1+\lambda_n)^2$ satisfies $c(\lambda_n, \delta t, L) > 0$. Then, we can derive the following inequality
\begin{equation}\label{ineq:quad_A}
    \frac{\alpha\lambda_{\min} + c(\lambda_{\max},\delta t, L)}{(1+\lambda_{\max})^2} \|u^{n+1} - u^n\|^2\leq A(u^n,u^{n-1}) - A(u^{n+1},u^n).
\end{equation}
Summing the inequality \eqref{ineq:quad_A} from $n=1$ to $\infty$, we finally obtain
\begin{equation*}
    \sum_{n=1}^{\infty}\|u^{n+1} - u^n\|^2\leq  \frac{(1+\lambda_{\max})^2}{\alpha\lambda_{\min} + c(\lambda_{\max},\delta t, L)}(A(u^1,u^0) - \liminf_{n\to \infty}A(u^{n+1},u^n))<\infty.
\end{equation*}
This completes the proof. 
\end{proof}
Now, let us turn to the global convergence of $\{ u^n\}_n$ by proving that the sequence  $\{ u^n\}_n$ has a finite length. 
\begin{theorem}\label{thm:globalconverg}
   Assume that $h$ and $f$ are Lipschitz continuous, $A(x,y) = E(x) + (\frac{L}{2}+\frac{1}{3\delta t})\|x-y\|^2$ is a KL function and  $ \lambda_{\max} <\sqrt{\frac{4}{3\delta tL} - \frac12}-1$. Let $\{u^n\}_n$ be a sequence generated by Algorithm \ref{alg:cap} with $\hat u^n=u^n$ and line search criterion \eqref{ineq:MS_LS_rule} for solving problem \eqref{eq:basefunctional}. Then the following statements hold.
\begin{itemize}
    \item[\rm (i)] $\displaystyle \lim_{n\to\infty}\|\nabla A(u^{n+1},u^{n})\|=0$.
    \item[\rm (ii)] The sequence $\{A(u^{n+1}, u^n)\}_n$ is monotone decreasing with a limit and the sequence $\{u^n\}_n$ is bounded. Specifically, there exists a constant $\zeta$ such that $\displaystyle \lim_{n\to \infty} A(u^{n+1},u^n)=\zeta$.
    \item[\rm (iii)] The sequence $\{u^n\}_n$ converges to a critical point of $E$; moreover,$\displaystyle\sum_{n=0}^{\infty}\|u^{n+1} - u^n\|<\infty$.
\end{itemize}
\end{theorem}
\begin{proof}
With directly calculating the gradient of $A(x,y)$ and notation $C_0 = {L} + \frac{2}{3\delta t}$, we have
\begin{equation}\label{eq:ls:grad_A}
    \nabla A(x,y)|_{(x,y)=(u^{n+1}, u^n)} = \left(
    \begin{array}{c}
        h(u^{n+1}) + f(u^{n+1}) + C_0(u^{n+1} - u^n)\\
    -C_0(u^{n+1} - u^n)
    \end{array}\right).
\end{equation}
Taking $\hat u^n = u^n$ in \eqref{eq:solve:y:pre:general}, we can obtain
\begin{equation*}
    f(u^{n+1}) = -\frac{2}{3\delta t}(3y^{n+1} - 4u^{n+1} + u^{n})  - \left(f(u^{n+1}) - f(u^{n})\right) -h(y^{n+1}) - M(y^{n+1} - u^{n+1}).
\end{equation*}
Substituting into \eqref{eq:ls:grad_A}, we find the first component of $ \nabla A(u^{n+1}, u^n)$ can be written as 
\begin{align*}
    &(h(u^{n+1}) - h(y^{n+1}))- (f(u^{n+1}) - f(u^{n}))-\frac{2}{3\delta t}(3y^{n+1} - 4u^{n+1}+ u^{n}) \\
    &- M(y^{n+1} - u^{n+1}) + C_0(u^{n+1} - u^n).
\end{align*}
Denote $L_h$ as the Lipschitz constant of $h$ and $\lambda_M$ as the largest eigenvalue of $M$ henceforth. We thus can estimate $\|\nabla A(u^{n+1},u^{n})\|$
\begin{align}
    & \|\nabla A(u^{n+1},u^{n})\|
    \leq  \|h(u^{n+1}) + f(u^{n+1}) + C_0(u^{n+1} - u^n)\| + C_0 \|u^{n+1} - u^n\| \notag  \\
    &\leq  \|h(u^{n+1}) - h(y^{n+1})- (f(u^{n+1}) - f(u^{n}))-\frac{2}{3\delta t}(3y^{n+1} - 4u^{n+1}+ u^{n}) - M(y^{n+1} - u^{n+1})\|  \notag \\
     &+ 2C_0 \|u^{n+1} - u^n\| \notag  \\
    &\leq \left(\frac{( \delta t L_h +2) }{\delta t (1+\lambda_{n+1})}+\frac{\lambda_M}{1+\lambda_{n+1}}\right) \|u^{n+2} -u^{n+1}\| + (L + \frac{2}{3\delta t} +2C_0) \|u^{n+1}-u^n\| \notag \\
    &\leq C_L(\|u^{n+2}-u^{n+1}\|+\|u^{n+1} - u^n\| ). \label{eq:esti:grad:A}
\end{align}
Here  $C_L = \max\left\{ L_h + \frac{2}{\delta t} +\lambda_M, L +\frac{2}{3\delta t} +2C_0\right\}$.
Since $\|u^{n+1}-u^n\| \rightarrow 0$ with Lemma \ref{lem:square:summable}, we conclude (i) by \eqref{eq:esti:grad:A}.

For (ii), by \eqref{ineq:quad_A}, we see that $A(u^n,u^{n+1})$ is decreasing with a lower bound.  The nonnegative sequence $\{A(u^n,u^{n+1})\}_n$ thus admits a limit $\zeta$. With the energy $E(x)$ being level-bounded on $x$, it can be checked that $A(x,y)$ is thus level-bounded on $(x,y)$. The boundedness of $u^n$ can be guaranteed since $A(u^n,u^{n+1})$ is bounded and level-bounded on $(u^n,u^{n+1})$.


With the assumption that $\psi(\cdot)$ is a concave function, we arrive at
\begin{align}\label{ineq:KL_of_psi}
    \notag&\left[\psi(A(u^{n+1}, u^n) - \zeta) - \psi(A(u^{n+2}, u^{n+1}) - \zeta)\right]\|\nabla A(u^{n+1},u^{n})\|\\ 
    & \geq  \psi'(A(u^{n+1}, u^{n}) - \zeta)\left[A(u^{n+1}, u^{n}) - A(u^{n+2}, u^{n+1})\right]\|\nabla A(u^{n+1},u^{n})\|\\
    \notag & \geq A(u^{n+1}, u^{n}) - A(u^{n+2}, u^{n+1})\geq\frac{\alpha\lambda_{\min} + c(\lambda_{\max}, \delta t, L)}{(1+\lambda_{\max})^2}\|u^{n+2}-u^{n+1}\|^2, 
\end{align}
where the second inequality employs $\psi'(A(u^{n+1}, u^{n}) - \zeta)\|\nabla A(u^{n+1},u^{n})\| \geq 1$ with the KL property of $A(x,y)$ and the third inequality follows \eqref{ineq:quad_A}.
For convenience, denoting
\begin{equation}\label{eq:notation:Phi}
\Psi(u^n, u^{n+1}, u^{n+2}, \zeta): = \psi(A(u^{n+1}, u^n) - \zeta) - \psi(A(u^{n+2}, u^{n+1}) - \zeta),
\end{equation}
we have
\begin{equation*}
    \|u^{n+2}-u^{n+1}\|^2\leq \Tilde{C}(\|u^{n+2} - u^{n+1}\| + \|u^{n+1} - u^{n}\|)\Psi(u^n, u^{n+1}, u^{n+2}, \zeta),
\end{equation*}
where $C = \frac{(1+\lambda_{\max})^2}{\alpha\lambda_{\min} + c(\lambda_{\max}, \delta t, L)}C_L$. Using the fact $a^2\leq cd \Rightarrow a\leq c + \frac{d}{4}$, we arrive at
\begin{align}\label{ineq:splitform}
   &\|u^{n+2} - u^{n+1}\| 
    \leq C\Psi(u^n, u^{n+1}, u^{n+2}, \zeta) + \frac14(\|u^{n+1} - u^n\| + \|u^{n+2} - u^{n+1}\|)   \notag\\ 
   & \frac12\|u^{n+2} - u^{n+1}\| 
    \leq C\Psi(u^n, u^{n+1}, u^{n+2}, \zeta) + \frac14(\|u^{n+1} - u^{n}\| - \|u^{n+2} - u^{n+1}\|).
\end{align}
Summing the above inequality from $n=T$ to $\infty$, remembering the definition of $\Psi$ in \eqref{eq:notation:Phi}, we finally obtain  
\begin{equation*}           
      \sum_{t=T}^{\infty}\|u^{n+1} - u^n\| \leq 2C\psi(A(u^{T}, u^{T-1}) - \zeta) + \frac12\|u^{T} - u^{T-1}\|<\infty.
\end{equation*}
Therefore, we can deduce that the sequence $\{u^n\}_n$ is a Cauchy sequence. Utilizing the Lipschitz continuity of functions $h$ and $f$, together with condition (i), we can establish that $\nabla E(u^{n+1})=h(u^{n+1}) + f(u^{n+1}) \rightarrow 0$, and the sequence $\{u^n\}_n$ converges to a stationary point of $E$. This, in turn, concludes our proof. 
\end{proof}
For the global convergence of the original second-order convex splitting methods \eqref{eq:bdf2:adam:orig} without preconditioners and line search, which have already been widely used in phase-field simulations or gradient flow  \cite{ARW1995, Feng2013, Li2017, Shen2010}, we have the following remark. 
\begin{remark}\label{rem:un:mul}
By setting $\lambda_n = 0$ and the proximal weight $M = 0$ as specified in Proposition \ref{prop:descent}, Lemmas \ref{lem:LSbase} and \ref{lem:square:summable}, and Theorem \ref{thm:globalconverg}, we can still obtain the global convergence for the original second-order convex splitting methods \eqref{eq:bdf2:adam:orig}. The global convergence can still be attained in cases where $\lambda_n = 0$ (implying preconditioning without line search) or $M = 0$ (indicating line search without preconditioners).
\end{remark}

The local convergence rate relies on the KL exponent of the Auxiliary function $A$. The following theorem and its proof are standard (see \cite[Theorem 2]{Attouch2009} or \cite[Lemma 1]{Artacho2018}), and we present them within the framework of our proposed algorithm.
\begin{theorem}[Local Convergence Rate]\label{them:local:rate}
Under the same assumptions as Theorem \ref{thm:globalconverg}, consider a sequence $\{u^n\}_n$ generated by Algorithm \ref{alg:cap} with $\hat u^n=u^n$ that converges to $u^*$. Assume that $A(x,y)$ is a $KL$ function with $\psi$ in the $KL$ inequality given as $\psi(s) = cs^{1-\theta}$, where $\theta\in[0,1)$ and $c>0$. The following statements are valid:
\begin{itemize}
\item[\emph{(i)}] For $\theta = 0$, there exists a positive integer $n_0$ such that $u^n$ remains constant for $n > n_0$.
\item[\emph{(ii)}] For $\theta\in (0, \frac12]$, there are positive constants $c_1$, $n_1$, and $\eta\in (0,1)$ such that $\|u^n - u^*\| < c_1\eta^n$ for $n > n_1$.
\item[\emph{(iii)}] For $\theta\in (\frac{1}{2},1)$, there exist positive constants $c_2$ and $n_2$ such that $\|u^n-u^*\| < c_2 n^{-\frac{1-\theta}{2
\theta-1}} $ for $n>n_2$.
\end{itemize}
\end{theorem}
\begin{proof}
First, we prove (i). If $\theta = 0$, we claim that there must exist $n_0 > 0$ such that $A(u^{n_0},u^{n_0-1}) = \zeta$. Suppose to the contrary that  $A(u^{n},u^{n-1}) > \zeta$ for all $n>0$. Since $\lim_{n\to \infty}u^n = u^*$ and the sequence $\{A(u^n,u^{n-1})\}_n$ is monotone decreasing and convergent to $\zeta$ by Theorem \ref{thm:globalconverg}(ii), we choose the concave function $\psi(s) = cs$ and the KL inequality \eqref{eq:kl:def} that for all sufficiently large $n$, $\|\nabla A(u^{n},u^{n-1})\|\geq c^{-1}$, which contradicts \rm{Theorem \ref{thm:globalconverg}\rm{(i)}}. Thus, there exist $n_0 > 0$ so that $A(u^{n_0},u^{n_0-1}) = \zeta$. Since the sequence $\{A(u^n,u^{n-1})\}_n$ is monotone decreasing and convergent to $\zeta$, it must hold that $A(u^{n_0+\Bar{n}},u^{n_0+\Bar{n}-1}) = \zeta$ for any $\Bar{n}>0$. Thus, we can conclude from \eqref{ineq:KL_of_psi} that $u^{n_0} = u^{n_0+\Bar{n}}$. This proves that if $\theta = 0$, there exists $n_0>0$ so that $u^n$ is constant for $n>n_0$.

Next, we consider the case that $\theta \in (0,1)$. Based on the proof above, we need only consider the case when $A(u^n, u^{n-1})>\zeta$ for all $n>0$.

Define $\Delta_n =A(u^n,u^{n-1}) - \zeta$ and $S_n = \sum_{i=n}^{\infty}\|u^{i+1} - u^{i}\|$, where $S_n$ is well-defined due to the Theorem \ref{thm:globalconverg}(iii). Then, from (\ref{ineq:splitform}), we have for any $n>N$ that 
\begin{align*}
    S_n = &2\sum_{i=n}^{\infty}\frac12\|u^{i+1} - u^i\|
    \leq 2\sum_{i=n}^{\infty}\left[\Tilde{C}\Psi(u^{i-1},u^i,u^{i+1},\zeta) + \frac14(\|u^{i} - u^{i-1}\| - \|u^{i+1} - u^{i}\|)\right]\\
    \leq&2\Tilde{C}\psi(A(u^{n}, u^{n-1}) - \zeta)  + \frac12\|u^{n} - u^{n-1}\|
    =2\Tilde{C}\psi(\Delta_n)  + \frac12(S_{n-1}-S_{n})\\
    \leq&2\Tilde{C}\psi(\Delta_n)  + \frac12(S_{n-1}-S_{n+1}).
\end{align*}
Set $ \psi(s) = cs^{1-\theta}$, for all sufficiently large $n$, $c(1-\theta)\Delta_n^{-\theta}\|\nabla A(u^{n+1},u^{n})\|\geq 1$, and we can transfer the estimate of $\|\nabla A(u^n,u^{n-1})\|$ to a new formulation
\begin{equation*}
   \|\nabla A(u^n,u^{n-1})\|\leq C_L(\|u^{n+1}-u^{n}\|+\|u^{n} - u^{n-1}\|) = C_L(S_{n-1}-S_{n+1}).
\end{equation*}
Due to the above two inequalities, we have
\begin{equation*}
    \Delta_n^{\theta}\leq (1-\theta) C_Lc(S_{n-1}-S_{n+1}).
\end{equation*}
Combining it with $S_n\leq2\tilde{C}\psi(\Delta_n)  + \frac12(S_{n-1}-S_{n+1})$, we have 
\begin{align}\label{eq:ineq_Sn}
    S_n\leq& 2c\tilde{C}(\Delta_n^{\theta})^{\frac{1-\theta}{\theta}} + \frac12(S_{n-1}-S_{n+1})\leq C_1(S_{n-1}-S_{n+1})^{\frac{1-\theta}{\theta}} + (S_{n-1}-S_{n+1})
\end{align}
where $C_1 = 2c\tilde{C}[(1-\theta)C_Lc]^{\frac{1-\theta}{\theta}}$. There are two cases: $\theta \in (0,\frac12]$ and $\theta \in (\frac12,1)$. Suppose the first that $\theta \in (0,\frac12]$, then $\frac{1-\theta}{\theta} \geq 1$. Since $\|u^{n+1} - u^n\| \to 0$ from Lemma \ref{lem:square:summable}, it leads to $S_{n-1} - S_{n+1} \to 0$. From these and \eqref{eq:ineq_Sn}, we can conclude that there exists $n_1>0$ such that for all $n\geq n_1$, we have
\begin{equation*}
    S_{n+1}\leq S_n\leq (C_1+1)(S_{n-1}-S_{n+1})
\end{equation*}
which implies that $S_{n+1} \leq \frac{C_1+1}{C_1+2}S_{n-1}$. Hence, for sufficiently large $n>n_1$
\begin{equation*}
    \|u^{n+1} - u^*\|\leq \sum_{i = n+1}^{\infty}\|u^{i+1} - u^i\| = S_{n+1}\leq S_{n_1-1}\eta^{n-n_1}, \quad  \eta := \sqrt{\frac{C_1 + 1}{C_1 + 2}}.
\end{equation*}
Finally, we consider the case that $\theta \in (\frac12,1)$, which imply $\frac{1-\theta}{\theta}\leq 1$. Combining this with the fact that $S_{n-2} - S_n \to 0$, we see that there exists $n_2 > 0$ such that for all $n\geq n_2$, we have
\begin{align*}
    S_{n+1}\leq S_n\leq& 2c{C}(\Delta_n^{\theta})^{\frac{1-\theta}{\theta}} + \frac{1}{2}(S_{n-1}-S_{n+1})
    =(C_1 + \frac{1}{2})(S_{n-1}-S_{n+1})^{\frac{1-\theta}{\theta}}.
\end{align*}
Raising both sides of the above inequality to the power of $\frac{\theta}{1-\theta}$, we can observe a new inequality that $S_{n+1}^{\frac{\theta}{1-\theta}} \leq C_2(S_{n-1}-S_{n+1})$ for $n \geq n_2$, where $C_2 = (C_1 + \frac12)^{\frac{\theta}{1-\theta}}$. Let us define the sequence $\Omega_n = S_{2n}$. For any $n \geq \left\lceil \frac{n_2}{2} \right\rceil$, with nearly the same arguments as in \cite[Page 15]{Attouch2009}, there exists a constant $C_3>0$ such that for sufficient large $n$ 
\begin{equation*}
    \Omega^{\frac{\theta}{1-\theta}}_{n} \leq C_2(\Omega_{n-1}- \Omega_{n}) \Rightarrow
    \Omega_{n} \leq C_3n^{-\frac{1-\theta}{2\theta -1}}.
\end{equation*}
This yields
\begin{align*}
    \|u^n - u^*\|
    \leq S_n\left\{
    \begin{array}{ll}
         = \Omega_{\frac{n}{2}}\leq 2^{\rho}C_3n^{-\rho} &\text{if $n$ is even}\\
         \leq \Omega_{\frac{n-1}{2}}\leq 2^{\rho}C_3(n-1)^{-\rho}\leq 4^{\rho}C_3n^{-\rho}   &\text{if $n$ is odd and $n\geq 2$}
    \end{array}
    \right.
\end{align*}
where $ \rho := \frac{1-\theta}{2\theta-1}$. The proof is completed. 
\end{proof}

\subsection{Convergence analysis of FlowBBApe in the nonsmooth setting}\label{sec:nonsmooth}
This section focuses on the global convergence of Algorithm \ref{alg:cap} for the case $\hat u^n = \tilde u^n = \frac{4}{3}u^n-\frac{1}{3}u^{n-1}$ with nonsmooth $H$ as introduced in Section \ref{sec:intro}. In this subsection, we employ the `NLS' criterion \eqref{ineq:ori_LS_rule} for the line search step.

Under the above assumption, the sequence $\{y^n\}_n$ obtained from the subproblem \eqref{eq:mini:dc:sub} obtained from the first-order optimality condition, 
\begin{equation}\label{eq:til:yn}
    \mathbf{0}\in\frac{2}{3\delta t}(3y^{n} - 4u^n + u^{n-1}) + M(y^{n} - \tilde u^n) + \partial H( y^{n}) + (2f(u^n) - f(u^{n-1})) .
\end{equation}
First, we derive the estimate of energy decay between $ y^{n-1} $ and $ y^n $. Using the auxiliary function $ \tilde{A} $, we establish that the sequence $ \{\|u^{n+1} - u^n\|_n\} $ is square summable. Consequently, we establish the same global convergence for the sequence $ \{u^n\}_n $ by using the KL property of the auxiliary function $ \tilde{A} $.


\begin{lemma}\label{lem:square:sum:til}
    Let $\{u^n\}_n$ be a sequence generated by Algorithm \ref{alg:cap} with $\hat u^n=\tilde u^n$. With the line search criterion \eqref{ineq:ori_LS_rule} and  $0\le\lambda_{n} <\min\left\{\sqrt{\frac{10-6\delta tL}{2+3\delta tL}}-1,\sqrt{5}-1\right\}$, the sequence $\{ \|u^{n+1} - u^n\|\}_n$ is square summable, i.e., 
    $\sum_{n = 1}^{\infty}\|u^{n+1} - u^n\|^2< \infty$.
\end{lemma}
\begin{proof}First, we discuss the relation between $E(y^n)$ and $E(y^{n-1})$ from the bound estimate of the term $E^n(y^n) - E(y^{n-1})$. It can be divided into three parts:
\begin{equation}\label{eq:divide_Ey}
    E^n(y^n) - E(y^{n-1}) = 
\underbrace{[E^n(y^n) - E^n(u^n)]}_{T_1} + 
\underbrace{[E^n(u^n) - E(u^n)]}_{T_2} + 
\underbrace{[E(u^n) - E(y^{n-1})]}_{T_3}.
\end{equation}

We obtain the upper bound of the first term $T_1$ from the Proposition \ref{prop:til:1} and the third term $T_3$ with the line search criterion \eqref{ineq:ori_LS_rule},
\begin{equation*}\label{ineq:esti_T_1_T_3}
    T_1 \le -\|y^n-u^n\|^2_{N_3} + \|u^n-u^{n-1}\|^2_{N_4},\quad T_3 \le  -\alpha\lambda_{n-1} \|y^{n-1}-u^{n-1}\|^2,
\end{equation*}
where $N_3 = \left( \frac{4}{3\delta t} - \frac{L}{2} \right)I + \frac56M$, $N_4 = \frac16M$. For the term $T_2$, with direct calculation, we have
\[
T_2=-\frac{1}{3\delta t}\|u^n-u^{n-1}\|^2 + \langle f(u^n)-f(u^{n-1}), u^n-u^{n-1} \rangle.
\]
Taking the formulation of $E^n$ into the equation \eqref{eq:divide_Ey}, we have
\begin{align}\label{ineq:relation_Ey}
& \underbrace{E(y^n)+ \frac{1}{\delta t} \|y^n - u^n\|^2 - \frac{1}{3\delta t} \|y^n - u^{n-1}\|^2 +\langle f(u^n)-f(u^{n-1}), y^n - u^{n-1}\rangle}_{E^n(y^n)} -E(y^{n-1})
  \notag\\
\le& -\frac{1}{3\delta t}\|u^n - u^{n-1}\|^2 + \langle f(u^n)- f(u^{n-1}),u^n-u^{n-1}\rangle - (\frac{4}{3\delta t} - \frac{L}{2})\|y^n - u^n\|^2\notag\\
&-\frac56\|y^n-u^n\|_M^2 + \frac16\|u^n - u^{n-1}\|^2_{M} - \alpha\lambda_{n-1}\|y^{n-1}-u^{n-1}\|^2.
\end{align}
With the inequality \eqref{ineq:relation_Ey}, we arrive at 
\begin{align*}
    &(\frac{7}{3\delta t}-\frac{L}{2}) \|y^n - u^n\|^2 + \frac56\|y^n-u^n\|_M^2 - \frac16\|u^n-u^{n-1}\|_M^2+\alpha\lambda_{n-1}\|y^{n-1}-u^{n-1}\|^2\\
    \le& E(y^{n-1}) - E(y^n) -\frac{1}{3\delta t}\|u^n - u^{n-1}\|^2 + \frac{1}{3\delta t} \|y^n - u^{n-1}\|^2 + \langle f(u^n)- f(u^{n-1}),u^n-y^n\rangle\\
    \le &E(y^{n-1}) - E(y^n) -\frac{1}{3\delta t}\|u^n - u^{n-1}\|^2 + \frac{1}{3\delta t} \|y^n - u^{n-1}\|^2 + \frac{L}{2}\|u^n - u^{n-1}\|^2 + \frac{L}{2}\|y^n-u^n\|^2\\
    \le&E(y^{n-1}) - E(y^n) -\frac{1}{3\delta t}\|u^n - u^{n-1}\|^2 + \frac{2}{3\delta t} (\|y^n -u^n\|^2+\|u^n- u^{n-1}\|^2) \\
    &+ \frac{L}{2}\|u^n - u^{n-1}\|^2 + \frac{L}{2}\|y^n-u^n\|^2
\end{align*}
where the first inequality combines the terms $\langle f(u^n) - f(u^{n-1}),\cdot\rangle$, the second inequality follows the Cauchy-Schwarz inequality and the Lipschitz continuity, and the third inequality follows the fact that $\|a+b\|^2\le2\|a\|^2 + 2\|b\|^2$.
We can now obtain the subsequent estimate,
\begin{align*}
    &(\frac{5}{3\delta t}-L) \|y^n - u^n\|^2 + \frac56\|y^n-u^n\|_M^2 - \frac16\|u^n-u^{n-1}\|_M^2+\alpha\lambda_{n-1}\|y^{n-1}-u^{n-1}\|^2\\
     \le&E(y^{n-1}) - E(y^n) + (\frac{1}{3\delta t}+\frac{L}{2})\|u^n- u^{n-1}\|^2. 
\end{align*}
With a new auxiliary function $\tilde A$, defined as
\begin{equation}\label{eq:auxi_func}
    \tilde A(x,y) := E(y) +(\frac{5}{3\delta t}-L)\|y-x\|^2+\frac56\|y-x\|_M^2.
\end{equation}
With $\|(u^n-u^{n-1}\|^2/(1+\lambda_{n-1})^2=\|y^{n-1}-u^{n-1}\|^2$ by the line search, we finally obtain the estimate,
\begin{align}\label{ineq:constraint_lambda}
    &\left(\frac{5-3\delta t L}{3\delta t(1+\lambda_{n-1})^2}-\frac{2+3\delta t L}{6\delta t}\right) \|u^n - u^{n-1}\|^2 + \left(\frac{5}{6(1+\lambda_{n-1})^2}-\frac16\right)\|u^n-u^{n-1}\|_M^2 \notag\\ &+\alpha\lambda_{n-1}\|y^{n-1}-u^{n-1}\|^2
     \le\tilde A(u^{n-1},y^{n-1}) - \tilde A(u^{n},y^{n}).
\end{align}
To ensure that the coefficients before $\|u^n - u^{n-1}\|^2$ and $\|u^n-u^{n-1}\|_M^2$ in \eqref{ineq:constraint_lambda} are positive, the constraint of the line search step $0\le\lambda_{n-1} <\min\{\sqrt{\frac{10-6\delta tL}{2+3\delta tL}}-1,\sqrt{5}-1\}$, there exists a constant $\tilde C_{\lambda}>0$ that
\begin{equation}\label{ineq:A_til_decay}
    \tilde C_{\lambda}\|u^{n} - u^{n-1}\|^2 \le \tilde A(u^{n-1},y^{n-1}) - \tilde A(u^{n},y^{n}).
\end{equation}
Summing the above inequality from $1$ to $\infty$, we have
\begin{equation*}
  \sum_{n=1}^{\infty}\tilde C_{\lambda}\|u^{n} - u^{n-1}\|^2 \leq \tilde A(u^0,y^0) - \liminf_{i\to\infty}  \left[E(y^i) +(\frac{5}{3\delta t}-L)\|y^i-x^i\|^2+\frac16\|y^i-x^i\|_M^2\right]<\infty.
\end{equation*}
This completes the proof.
\end{proof}
\begin{remark}
    In inequality \eqref{ineq:constraint_lambda}, the estimate is valid only when the coefficients are positive. We consider the scenario where there is no energy decay along the search direction $d^n$ at the $n$-th step. In this case, the search step size satisfies $\lambda_n = 0$, so we exclude the term involving the coefficient $\alpha\lambda_{n-1}$. Moreover, the upper bound is well defined, since ${10-6\delta tL}>{2+3\delta tL}$ holds under the initial constraint $\delta t < \frac{2}{3L}$.
\end{remark}
\begin{theorem}\label{thm:globalconverg:M}
    Assume that $\tilde A(x,y) = E(y) +(\frac{5}{3\delta t}-L)\|y-x\|^2+\frac56\|y-x\|_M^2$ is a KL function. Let $\{u^n\}_n$ be a sequence generated by Algorithm \ref{alg:cap} with $y^n$ updated by \eqref{eq:til:yn}. Under the same assumptions of Lemma \ref{lem:square:sum:til}, the following statements hold:
\begin{itemize}
    \item[\rm (i)] $\displaystyle \lim_{n\to\infty} \dist(\mathbf{0}, \partial \tilde A(u^{n},y^{n}))=0$.
    \item[\rm (ii)] The sequence $\{\tilde A(u^{n},y^{n})\}_n$ is monotonically decreasing with a limit and the sequence $\{u^n\}_n$ is bounded. Specifically, there exists a constant $\tilde \zeta$ such that $\displaystyle \lim_{n\to \infty} \tilde A(u^{n},y^{n})=\tilde \zeta$.
    \item[\rm (iii)] The sequence $\{u^n\}_n$ converges to a stationary point of $E$; moreover, $\displaystyle\sum_{n=0}^{\infty}\|u^{n+1} - u^n\|<\infty$.
\end{itemize}
\end{theorem}

\begin{proof}
With directly calculating the subgradient of $\tilde A(x,y)$ at $(u^{n},y^{n})$, we have
\begin{equation*}
   \partial \tilde A(u^{n},y^{n}) =\left( 
    \begin{array}{c}
         D_1(u^{n} - y^{n} ) + \frac{5}{3}M(u^{n}-y^{n} )\\
         \partial H(y^{n}) + f(y^{n})+D_1(y^{n} - u^{n} ) + \frac{5}{3}M(y^{n}-u^{n} )
    \end{array}\right)
\end{equation*}
where $D_1 = \frac{10}{3\delta t} - 2L $. Based on the first-order optimality condition \eqref{eq:til:yn}, there exists $v_n\in \partial H(y^n)$ satisfying
\begin{equation}\label{eq:first_order_opt_condition}
        0 = \frac{2}{3\delta t}(3y^{n} - 4u^{n} + u^{n-1}) +v^n + \left(2f(u^{n}) - f(u^{n-1})\right)  + M(y^{n} - u^{n}).
\end{equation}
Utilizing the maximum eigenvalue $\lambda_M$ of the positive-semidefinite matrix $M$ along with the Lipschitz continuity of the function $f$, we proceed to estimate the distance as follows:
\begin{align}\label{eq:esti:subgrad:A:til}
    &\text{dist}(\mathbf{0},\partial \tilde A(u^n,y^n))\notag\\
    \le &2(D_1+\frac13\lambda_M)\|y^{n}-u^{n}\| + (\frac{2}{\delta t}+L)\|y^n-u^n\| +(\frac{2}{3\delta t}+L) \|u^n - u^{n-1}\|+ \lambda_M\|y^n - u^n\|\notag\\
    \le&E_1\|u^{n+1} - u^n\| + E_2\|u^n - u^{n-1}\|\le \tilde C_L(\|u^{n+1} - u^n\| +\|u^n - u^{n-1}\|)
\end{align}
where $C_1 = \frac{2+\delta t (2D_1+L)}{\delta t(1+\lambda_n)}+\frac{5+2\lambda_n}{3(1+\lambda_n)}\lambda_M$, $C_2 = \frac{2}{3\delta t}+L$ and $\tilde C = \max(C_1,C_2)$. Since $\|u^{n+1}-u^n\| \rightarrow 0$ with Lemma \ref{lem:square:sum:til}, we conclude (i) by \eqref{eq:esti:subgrad:A:til}.

For (ii), according to \eqref{ineq:A_til_decay}, we observe that $\tilde A(u^n,y^n)$ is decreasing with a lower bound. The sequence $\{\tilde A(u^n,y^n)\}_n$ thus admits a limit $\tilde \zeta$. The boundedness of $u^n$ can be guaranteed since $\tilde A(u^n,y^n)$ is bounded and level-bounded on $(u^n,y^n)$.

With the assumption that $\psi(\cdot)$ is a concave function, we arrive at
\begin{align}
    \notag&\left[\psi(\tilde A(u^n,y^n) - \tilde \zeta) - \psi(\tilde A(u^{n+1},y^{n+1}) - \tilde \zeta)\right]\text{dist}(\mathbf{0}, \partial \tilde A(u^{n},y^{n}))\\ 
    & \geq  \psi'(\tilde A(u^n,y^n) -  \tilde\zeta)\left[\tilde A(u^n,y^n) - \tilde A(u^{n+1},y^{n+1})\right]\text{dist}(\mathbf{0}, \partial \tilde A(u^{n},y^{n}))\\
    \notag & \geq \tilde A(u^n,y^n) - \tilde A(u^{n+1},y^{n+1})\ge \tilde C_{\lambda}\|u^{n+1} - u^{n}\|^2
\end{align}
where the second inequality uses the KL property of $\tilde A(x,y)$ and the third inequality follows \eqref{ineq:A_til_decay}. Employing the same methodology as described in Theorem \ref{thm:globalconverg}, we ultimately derive the following estimate,
\begin{equation*}           
      \sum_{t=T}^{\infty}\|u^{n+1} - u^{n}\| \leq 2\tilde{C}\psi(\tilde A( u^{T},y^{T}) - \zeta) + \frac12\|u^{T} - u^{T-1}\| <\infty,
\end{equation*}
where $\tilde C = \frac{\tilde C_L}{\tilde C_{\lambda}}$. This estimate implies that the sequence $\{u^n\}_n$ is convergent. This completes the proof.
\end{proof}

Let us conclude this section with a remark regarding the convergence of two special cases and the local convergence rate.
\begin{remark}
    For the case $\hat u = u^n$, a similar conclusion holds with the  auxiliary $\hat A(x,y) = E(y) + (\frac{5}{3\delta t}-\frac{L}{2})\|x-y\|^2$ and the step constraint $0\le\lambda_{n}<\sqrt{\frac{5-3\delta t L}{2+3\delta tL}}-1$. For the case $M=0$, there is no distinction between $\hat u^n= u^n$ or $\hat u^n = \tilde u^n$. Moreover, the local convergence rate analysis relies on the inequality \eqref{ineq:A_til_decay}. The auxiliary function $\tilde A$ retains the same properties as the earlier auxiliary function $A$, and the analogous results of Theorem \ref{them:local:rate} hold in this subsection.
\end{remark}

\section{Preconditioners and numerical experiments} \label{sec:pre}
\subsection{The preconditioning techniques along with the order of preconditioned scheme}\label{sec:pre:sub}
We now focus on the preconditioning technique applied to the case where $h(u) = Au - b_0$, where $A$ is a linear, positive-semidefinite, bounded operator and $b_0 \in X$ is known. The core idea is to employ classical preconditioning methods, including classical symmetric Gauss-Seidel, Jacobi, and Richardson preconditioners, to effectively deal with large-scale linear systems. Notably, we prove that any finite and feasible preconditioned iterations can ensure the global convergence of the whole nonlinear second-order convex splitting method. For a more comprehensive discussion of preconditioning techniques for nonlinear convex problems, see \cite{BSCC}; for further insights on DCA, see \cite{DS, Shensun2023}. 
To illustrate the fundamental principle of preconditioning techniques, we present the classical preconditioning technique via the following proposition.

\begin{proposition}\label{prop:proximal_to_pre}
If $h(u)=Au -b_0$, the equation \eqref{eq:solve:y:pre:general} can be rewritten as the preconditioned iteration 
\begin{equation}
    y^n =  \hat{u}^n + \mathbb{M}^{-1}(b^n - T\hat{u}^n),
\end{equation}
where 
\[
T = \frac{2}{\delta t}I + A, \quad \mathbb{M}= T +M, \quad b^n = b_0 + \frac{2}{3\delta t}(4u^n-u^{n-1}) -(2f(u^n) -f(u^{n-1})).
\]
This corresponds to one-step classical preconditioned iteration for solving \eqref{eq:bdf2:adam:orig}, i.e., $Ty^n = b^n$ with initial value $\hat u^n$.
\end{proposition}
\begin{proof}
    From \eqref{eq:solve:y:pre:general}, we have 
    \begin{equation}
            (\frac{2}{\delta t}I + A)y^n + My^n = M \hat{u}^n + b_0 + \frac{2}{3\delta t}(4u^n-u^{n-1}) -(2f(u^n) -f(u^{n-1})).
    \end{equation}
Using the definitions of $T$, $\mathbb{M}$, and $b^n$ as in the proposition, we obtain
\begin{align*}
Ty^n + My^n = M\hat{u}^n  + b^n 
&  \Leftrightarrow (T+M)y^n  = (T+M)\hat{u}^n  + b^n - T \hat{u}^n \\
&  \Leftrightarrow  y^n =  \hat{u}^n + \mathbb{M}^{-1}(b^n - T\hat{u}^n)
\end{align*}
which completes the proof.
\end{proof}

Note that classical preconditioners such as symmetric Gauss-Seidel (SGS) preconditioners do not require the explicit definition of $M$ in \eqref{eq:mini:dc:sub} \cite{BSCC}. The positive-semidefiniteness of $M$ is inherently fulfilled through SGS iteration in solving the linear equation $Ty = b^n$. 

Supposing $h(u) = Au-b_0$ as in Proposition \ref{prop:proximal_to_pre}, the discrete matrix of the linear operator $T=\frac{2}{\delta t}I + A$ is $\mathbb{D}-\mathbb{E}-\mathbb{E}^T$ where $\bD$ is the diagonal part and $-\bE$ represents the strictly lower triangular part. For a symmetric Gauss-Seidel (SGS) preconditioner, it is known that the precondition $\mathbb{M}= (\bD-\bE)\bD^{-1}(\bD-\bE^T)$ \cite[Section 4.1.2]{SA}. We thus obtain the weight in \eqref{eq:mini:dc:sub} $M=\mathbb{M}-T=\bE\bD^{-1}\bE^T$ is only positive-semidefinite.

Regarding the order of the term $M(u^{n+1} - u^{n})$ in \eqref{eq:solve:y:pre:general}, the following proposition shows that the SGS preconditioners remain of order $O((\delta t)^2)$. Thus, the preconditioned equation \eqref{eq:solve:y:pre:general} retains the second-order accuracy of the original second-order BDF and Adams-Bashforth scheme \eqref{eq:bdf2:adam:orig}. 
\begin{proposition}\label{prop:second_order:pre}
    For the SGS preconditioner and small $\delta t$, assume $\|d_n\|=\|y^n-u^n\| \leq C\delta t $ for all $n$. Then the order of the term $\|M(u^{n+1} - \hat u^n)\|$ in \eqref{eq:solve:y:pre:general} is $O((\delta t)^2)$.
\end{proposition}
\begin{proof}
Since $u^{n+1}=u^n + (1+\lambda_n)d_n$ and $u^{n-1}= u^{n} - (1+\lambda_{n-1})d_{n-1}$, the assumption on $\|d_n\|$ implies
\begin{equation}\label{eq:error:differ:hatu}
\|u^{n+1}-u^{n}\| = O(\delta t), \quad \|u^{n+1}-\tilde u^{n}\| = O(\delta t).
\end{equation}
Note that $T=\frac{2}{\delta t}I + A$, so $\|\bD\| = O(\frac{1}{\delta t})$ and $\|\bE\|= \|\bE^T\|= O(1)$. We compute
\begin{align}
    M &= \mathbb{M} - T = (\bD- \bE)\bD^{-1}(\bD-\bE^T) - \bD + \bE + \bE^T \notag \\
    & =\bD - \bE - \bE^T + \bE\bD^{-1}\bE^T - \bD + \bE + \bE^T \notag \\
    & = \bE\bD^{-1}\bE^T  \Rightarrow \| \bE\bD^{-1}\bE^T\| =O(\delta t). \label{eq:pre:order}
\end{align}
We derive this proposition by combining the equations \eqref{eq:error:differ:hatu} and \eqref{eq:pre:order}. 
\end{proof}

SGS preconditioners have already been used for large-scale linear subproblems in convex, nonsmooth, and nonlinear optimization \cite{BSCC}. General Jacobi preconditioners incorporating $\mathbb{M}=c \bD$ subject to the condition $\mathbb{M} \succeq T$ prove to be highly beneficial for parallel computing \cite{Shensun2023}. However, the task of selecting an appropriate value for $c$ to ensure that $M = \mathbb{M}-T$ is of the order $O(\delta t)$ poses a significant challenge. By opting for $\mathbb{M}=\frac{2}{\delta t}I + \tilde c \text{Diag}(A)$ under the condition $\tilde c \text{Diag}(A) \succeq A$, it becomes evident that $\mathbb{M} \succeq T$ and is indeed feasible. It can be verified that $\|M\| = \|\mathbb{M}-T\| = O(1)$ in this case. Indeed, for general Jacobi preconditioners, the modified equation \eqref{eq:solve:y:pre:general} fails to maintain the second-order characteristic of the original second-order BDF and Adams-Bashforth scheme that the SGS preconditioners in Proposition \ref{prop:second_order:pre} can achieve. However, the global convergence can be guaranteed for Jacobi preconditioners throughout this study.

\subsection{Numerical experiments}
In this section, we present two numerical experiments that validate the effectiveness of our algorithm in solving two typical nonconvex models. Numerical experiments for the least-squares problem with the SCAD regularizer are performed on a workstation equipped with an Intel(R) Xeon(R) CPU E5-2699A v4 @ 2.20 GHz. For the segmentation problem based on the nonlocal Ginzburg–Landau model, the code is implemented on a computer with an NVIDIA RTX 3080Ti (16 GB) GPU. Further implementation details and experimental results are detailed in the following subsections. For the KL properties and corresponding global convergence of these two models, see Remark \ref{rem:KL:conver:two:models} at the end of this section.
\subsubsection{Least squares problems with SCAD and modified SCAD regularizer}
 We consider the smoothly clipped absolute deviation (SCAD) regularization and a modified SCAD regularization. For the KL property of the auxiliary function and the SCAD regularization applied on the least squares problem, see \cite[Corollary 4.1]{Liu2019}. The DC decomposition of the original SCAD can be expressed as follows (see \cite[Section 6.1]{APX} or \cite{Wen2018})
\begin{equation}\label{eq:scad_regularization}
    P(u) = \mu\sum_{i=1}^{k}\int_0^{|u_i|}\min\left\{1,\frac{[\theta\mu-x]_+}{(\theta-1)\mu}\right\}dx=\mu\|u\|_1-\mu\sum_{i=1}^{k}\int_0^{|u_i|}\frac{[\min\{\theta\mu,x\}-\mu]_+}{(\theta-1)\mu}dx
\end{equation}
where $\theta > 1$ is a constant, $\mu>0$ represents the regularization parameter and $[x]_{+} = \max\{0,x\}$. 
 We denote the second term of SCAD, i.e., the integral term, as $P_2(u) = \mu \|u\|_1 -P(u) = \sum_{i=1}^k p_2(u_i)$. This term is continuously differentiable with the explicit expression and gradient
 \begin{equation*}
     p_2(u_i) = \left\{
        \begin{array}{ll}
         0& \text{if}  \ \  \ |u_i|\leq\mu \\
         \frac{(|u_i| - \mu)^2}{2(\theta - 1)}& \text{if} \ \  \ \mu<|u_i|<\theta\mu \\
         \mu|u_i| - \frac{\mu^2(\theta+1)}{2}&\text{if}\ \ \ |u_i|\geq\theta\mu 
    \end{array},
        \right. \ \ \nabla_i P_2(u_i)= \text{sign}(u_i) \dfrac{[\min\{\theta \mu,|u_i|\} -\mu]_{+}} {(\theta-1)}.
\end{equation*}
To ensure differentiability in the DC decomposition, we modify the $l_1$ norm using the Huber function $\mathcal{H}_{\alpha}(u):= \sum_{i}\mathcal{H}(u_i,\alpha)$, denoted as
\begin{equation*}
    \mathcal{H}(u_i,\alpha) = \left\{
    \begin{array}{ll}
         \frac{|u_i|^2}{2\alpha}& \text{if}  \ \  \ |u_i|\leq\alpha \\
         |u_i|-\frac{\alpha}{2}& \text{if} \ \  \ |u_i|>\alpha
    \end{array}.
    \right. 
\end{equation*}
Here, the scale $\alpha$ is a shape parameter for controlling the level of robustness, and we set $\alpha$ to be equal to $0.5\mu$ in the following experiment.
We employ two types of SCAD regularization in the least squares problem:
\begin{equation}\label{eq:lsp}
    \min_{u\in \mathbb{R}^k}E(u) = \frac12\|Au-b\|^2 + \mu P_1(u) - P_2(u),
\end{equation}
where $A\in \mathbb{R}^{m\times k}$, $b\in \mathbb{R}^{m}$. It can be checked that the modified SCAD regularization $P_{M}(u): =  \mu\mathcal{H}_{\alpha}(u) - P_2(u) := \sum_{i=1}^k p_{M}(u_i)$, where
   $p_{M}(u_i)$ can be expressed as  piecewise formulation:
        \begin{equation}\label{eq:huber-scad}
        \frac{p_{M}(u_i)}{\mu} = \left\{
        \begin{array}{ll}
       |u_i|^2/(2\alpha)& \text{if}  \ \  \ |u_i|\leq\alpha \\
         |u_i|-{\alpha}/{2}& \text{if}  \ \  \ \alpha < |u_i|\leq\mu \\
         |u_i|-{\alpha}/{2} -(|u_i| - \mu)^2/(2(\theta - 1)\mu)& \text{if} \ \  \ \mu<|u_i|<\theta\mu \\
         (\mu(\theta+1) - \alpha)/{2} &\text{if}\ \ \ |u_i|\geq\theta\mu 
    \end{array}.
        \right.
    \end{equation}
    
To construct the test problem, we generate a randomly normalized matrix $A$ (unit-norm columns) with dimensions $m \times k$. We then construct an index set $T \subset \{1, 2, 3, ..., k\}$ of size $s$ and build an $s$-sparese vector $y$ supported on $T$. Finally, we define the vector $b = Ay + 0.01\cdot\hat{k}$, where $\hat{k}$ consists of i.i.d. Gaussian random entries. All algorithms are initialized from the origin and terminate when
$$
\frac{\|u^n - u^{n-1}\|}{\max\{1,\|u^n\|\}}<10^{-12}.
$$

Evaluating various algorithms involves recording the iteration count (iter) and CPU time. All algorithms show slight variations in the final function value (fval) and the sparsity of the solution (sparsity). We use two single columns to record the `fval' and `sparsity'. In our numerical experiments, we track the sparsity by recording the count of values exceeding $10^{-5}$ in the solution. To assess their performance, we focus on problem sizes represented by $(m,k,s) = (720i,2560i,80i)$, where $i = 1,2,3,...,10$. The computational results are presented in Table \ref{table:lsp}, which correspond to the problem \eqref{eq:lsp} with $\mu = 3\times10^{-2}$ and $\theta = 10$. For each problem size, 5 instances are randomly generated, and the reported value for each line in Table \ref{table:lsp} represents the average across all the cases. We compare six algorithms to solve \eqref{eq:lsp}: our algorithms FlowBAp, FlowBApe, FlowBBAp, and FlowBBApe, along with the DC algorithm and the BDCA-Backtracking algorithm as presented in \cite{Artacho2018}. In alignment with our analysis, we employ the `MPCLS' line search criterion for the smooth case ($P_1 = \mathcal{H}_{\alpha}$). The `NLS' criterion is utilized for the nonsmooth case ($P_1 = \|\cdot\|_1$). We discuss the implementation details of these algorithms below.
\begin{itemize}
    \item \textbf{FlowBBAp} (\textbf{FlowBBApe}): This algorithm is represented by Algorithm \ref{alg:cap} with $\hat u^n= u^n$ (FlowBBAp) or $\hat u^n=\tilde u^n$ (FlowBBApe), where $H(u) = \mu \mathcal{H}_{\alpha}(u) + \frac12\|Au-b\|^2$ (or nonsmooth function $H(u) = \mu\|u\|_1 + \frac12\|Au-b\|^2$) and $F(u) = - P_2(u)$. The Lipschitz constant $L$ of $\nabla F$ is $1/(\theta-1)$ \cite[Example 4.3]{Wen2018} with $L=1/9$. The parameter values for this algorithm are chosen as follows: $\delta t < 2/(3L)$ (e.g., $\delta t= 6-10^{-15}$), $\bar \lambda_{\text{max}} = 5$, $\Bar{\lambda} = 0.618\bar \lambda_{\text{max}}$, $\alpha = 0.2$, $\beta = 0.8$, and $M = \lambda_{A^TA} I- A^TA$, where $\lambda_{A^TA}$ represents the largest eigenvalue of the matrix $A^TA$.
    \item \textbf{FlowBAp} (\textbf{FlowBApe}): This is a basic version of the algorithm {FlowBBAp} (or {FlowBBApe}) without the boosted method. The parameters of this algorithm remain consistent with those of the algorithm {FlowBBAp} (or {FlowBBApe}).
    \item \textbf{BDCA}: This is a boosted algorithm based on the DC algorithm. It further accelerates the fundamental algorithm, DCA, by employing the line search technique, as described in \cite[Algorithm 2]{Artacho2018}. To simplify computations, we employ distinct convex splitting to sidestep equation solving. The convex splitting can be defined as follows:
    \begin{equation} \label{eq:another:split}
        E(u) =\mu \mathcal{H}_{\alpha}(u) + \frac{\lambda_{A^TA}}{2}\|u\|^2 - \left(\frac{\lambda_{A^TA}}{2}\|u\|^2 +  P_2(u)  - \frac12\|Au-b\|^2\right)
    \end{equation}
    where $\lambda_{A^TA}$ represents the maximum eigenvalue of the matrix $A^TA$. The parameter of the line search part is the same as the algorithm FlowBBAp.
    \item \textbf{DCA}: This is the classical DC algorithm, which is a special version of the algorithm BDCA-Backtracking without the line search \cite[Algorithm 1]{Artacho2018}. 
\end{itemize}

From Table \ref{table:lsp}, it is evident that our proposed FlowBBApe algorithm always outperforms other algorithms. It exhibits clear advantages in iteration count and CPU time, especially for large-scale problems. Preconditioning alleviates step-size restrictions and enhances flexibility. Furthermore, with  $M = \lambda_{A^TA} I- A^TA$, one can choose $F(u) = -P_2(u)$, $\delta t < 2/(3L)=6$ with $L$ being the Lipschitz constant of $P_2$. However, for computing explicit resolvents, one has to choose $F(u) = \frac{1}{2}\|Au-b\|^2 - P_2(u)$ with $\delta t < 2/(3\lambda_{A^TA}) < 1/12$ since $\lambda_{A^TA}>8$, where the step size is much smaller.

\begin{sidewaystable}[htbp!]
\caption{Solving \eqref{eq:lsp} on random instances(\protect\textcircled{1}DCA \protect\textcircled{2}\protect\textcircled{7}BDCA \protect\textcircled{3}FlowBAp \protect\textcircled{4}FlowBApe \protect\textcircled{5}FlowBBAp \protect\textcircled{6}\protect\textcircled{8}FlowBBApe)\label{table:lsp}}
\begin{tabular}{cc|rrrrrr|rr|rrrrrr|rr|c|c}
\hline
\multicolumn{2}{c|}{Size} & \multicolumn{8}{c|}{Iter} & \multicolumn{8}{c|}{CPU time (s)} &Fval  &Sparsity \\ \hline
 && \multicolumn{6}{c|}{Smooth} & \multicolumn{2}{c|}{Nonsmooth} & \multicolumn{6}{c|}{Smooth} & \multicolumn{2}{c|}{Nonsmooth}  & \multicolumn{2}{c}{\multirow{2}{*}{Aver.}}\\ 
i&s &\textcircled{1}  & \textcircled{2} & \textcircled{3}  & \textcircled{4} & \textcircled{5}  & \textcircled{6} & \textcircled{7}  & \textcircled{8}  &\textcircled{1}  & \textcircled{2} & \textcircled{3}  & \textcircled{4} & \textcircled{5}  & \textcircled{6} & \textcircled{7}  & \textcircled{8}\\ \hline
1&80 &  590 &   313 &   603 &   406 &  259 &  \textbf{149} &  319 &  \textbf{187} &  \textbf{ 0.1} &   0.2 &   \textbf{0.1} &   \textbf{0.1} &   0.3 &   0.2 &   0.2 &   \textbf{0.1}& 0.389  &81.8\\
2&160&  637 &   340 &   651 &   439 &  270 &  \textbf{167} &  335 &  \textbf{190} &   \textbf{0.2} &   0.5 &   0.3 &   \textbf{0.2} &   0.5 &   0.3 &   0.5 &  \textbf{0.3} & 0.785  &169.2\\
3&240&  636 &   331 &   650 &   439 &  271 &  \textbf{158} &  335 &  \textbf{170} &   \textbf{2.4} &   4.5 &   3.9 &   2.7 &   4.9 &   3.0 &   4.6 &   \textbf{2.7}& 1.181  &251\\
4&320&  637 &   348 &   652 &   440 &  268 &  \textbf{160} &  361 &  \textbf{178} &   6.8 &  12.7 &   7.8 &   \textbf{5.0} &   9.1 &   5.4 &   9.5 &   \textbf{5.6} & 1.582  &334\\
5&400&  636 &   335 &   649 &   438 &  260 &  \textbf{168} &  324 &  \textbf{185} &  12.2 &  21.3 &  13.9 &   \textbf{8.2}&  13.7 &   9.4 &  13.9 &   \textbf{9.2} & 1.949   &427\\
6&480&  619 &   336 &   633 &   427 &  263 &  \textbf{168} &  345 &  \textbf{181} &  12.3 &  23.5 &  17.6 &  \textbf{11.6} &  21.2 &  13.9 &  23.5 &  \textbf{15.0} & 2.344   &498.8\\
7&560&  639 &   345 &   654 &   441 &  266 &  \textbf{162} &  342 &  \textbf{189} &  16.0 &  30.8 &  23.6 &  \textbf{16.0} &  27.1 &  16.6 &  27.3 &  \textbf{17.9} & 2.745  &482.8\\
8&640&  640 &   344 &   654 &   441 &  275 &  \textbf{165} &  338 &  \textbf{180} &  21.5 &  40.5 &  31.8 &  \textbf{20.9} &  38.6 &  22.8 &  37.2 &  \textbf{23.0}  & 3.151 &678\\
9&720&  631 &   356 &   645 &   436 &  267 &  \textbf{166} &  358 &  \textbf{179} &  28.4 &  56.1 &  42.0 &  \textbf{26.2} &  50.1 &  31.5 &  54.2 &  \textbf{32.2} & 3.525   &752.6\\
10&800&  631 &   339 &   645 &   435 &  265 &  \textbf{166} &  355 &  \textbf{186} &  31.9 &  58.4 &  44.3 & \textbf{ 27.1} &  44.8 &  28.5 &  53.6 &  \textbf{32.0} & 3.917  &831.4\\\hline
\end{tabular}
\end{sidewaystable}

\begin{figure}[htbp]
    \centering
    \begin{minipage}[b]{0.48\textwidth}
        \centering
        \begin{overpic}[width=\linewidth]{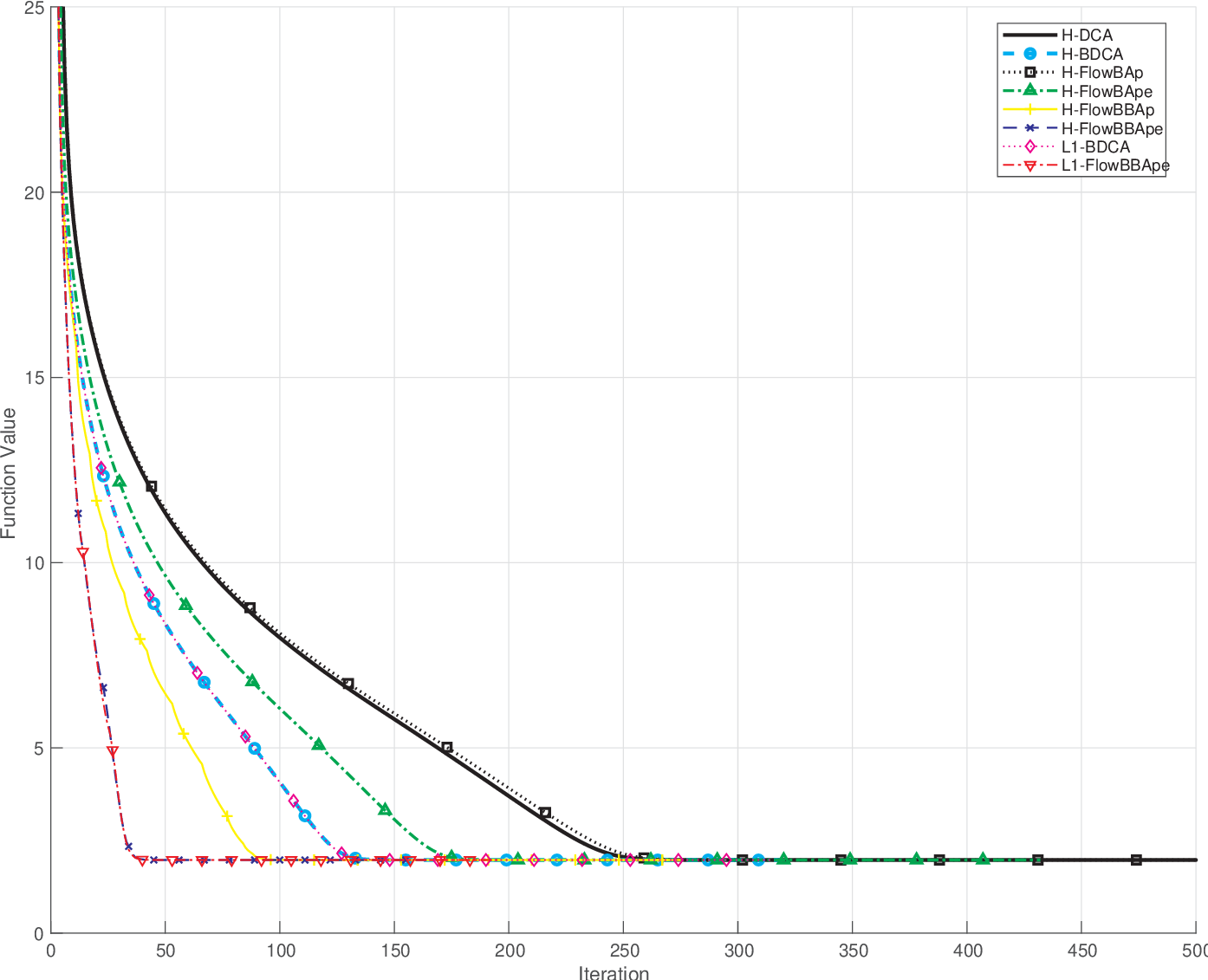}
            \put(40,17){\includegraphics[width=0.6\linewidth]{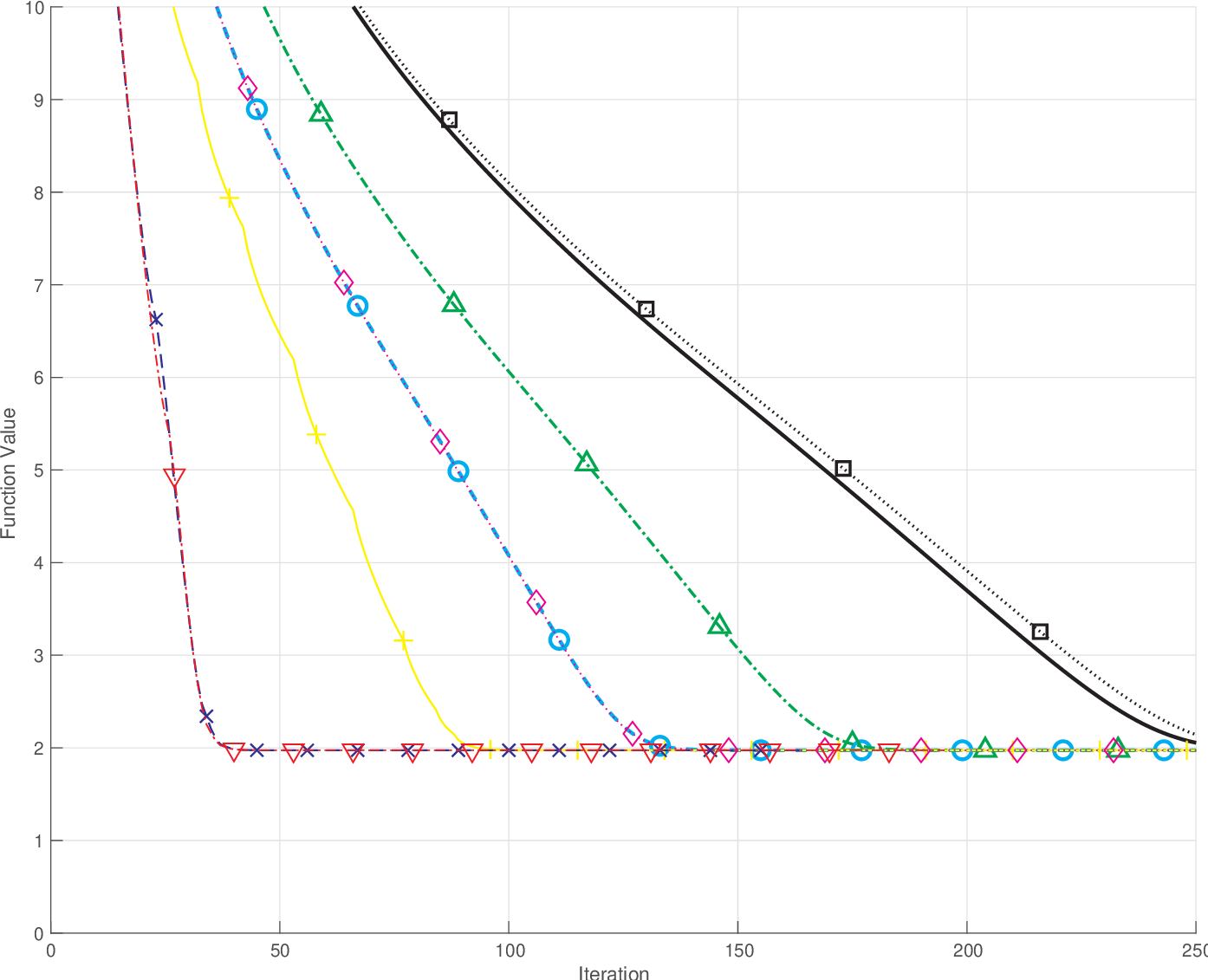}}
        \end{overpic}
        \caption{\footnotesize Convergence performance (`L1-' means $P_1(u) = \mu\|u\|_1$ and `H-' means $P_1(u) = \mu\mathcal{H}_{\alpha}(u)$)}
        \label{fig:cong_curves}
    \end{minipage}
    \hfill
    \begin{minipage}[b]{0.48\textwidth}
        \centering
        \begin{overpic}[width=\linewidth]{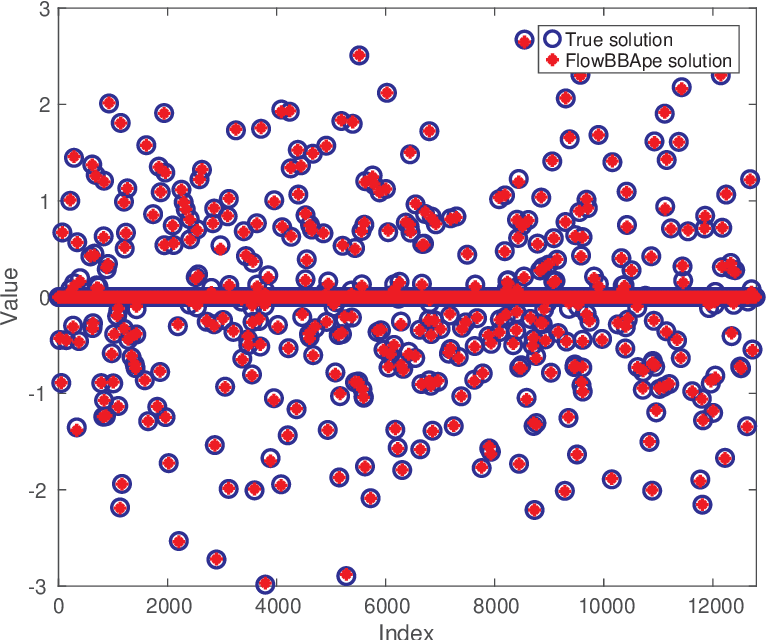}
        \end{overpic}
        \caption{\footnotesize The true solution and the solution obtained by \eqref{eq:lsp} with $P_1(x)=\mu\|x\|_1$ and $\mu=0.03$}
        \label{fig:sol_sp}
    \end{minipage}
\end{figure}

In Table \ref{table:lsp}, we present two additional experiments involving the least squares problem with the original SCAD regularization \eqref{eq:scad_regularization}. The comparison is made between our algorithm, FlowBBApe, and BDCA, which is an existing DCA algorithm utilizing line search \cite{Artacho2018}. These additional experiments further demonstrate the superior effectiveness of our algorithm in dealing with nonsmooth objective functions. In Figure \ref{fig:cong_curves}, we present the performance comparison of the various algorithms discussed earlier. Our algorithm, FlowBBApe, exhibits a significant decrease in the objective function value under both conditions: $P_1(u) = \mu\|u\|_1$ and $P_1(u) = \mu\mathcal{H}_{\alpha}(u)$. Figure \ref{fig:sol_sp} illustrates the disparity between the true solution $y$ and the solution derived from solving \eqref{eq:lsp} with the parameter $i=5$. The majority of the values in the FlowBBApe solution closely match the true solution, indicating its effectiveness as an approximation.
\subsubsection{Graph Ginzburg-Landau model}
In this subsection, we focus on a segmentation problem that utilizes the graph Ginzburg-Landau modeling. 
First, we introduce the graph structure, which is important for the construction of the graph  
Laplacian operator for the model. It contains the relationships among the image's pixels or  
data. We can regard the image or data as an undirected weighted graph $ G = (V, E) $, where $ V $ is  
the set of nodes corresponding to the image's pixels or the data points, and $ E \subset V \times V $ is the  
set of edges. The edge weighting is a function $ w: E \to \mathbb{R} $ corresponding to the weights  
set $ \{w_{ij}\} $ with $ w_{i,j} = w(e_{i,j}) \geq 0 $ and the edge $ e_{i,j} \in E $. For more details about this model, please refer to \cite{Shensun2023, BF}. Unlike the conventional phase-field model, we incorporate a prior term to complete the semi-supervised assignment. The formulation for this problem is expressed as follows \cite{BF}:
\begin{equation}\label{eq:nonlocalGL}
    \min_{u\in \mathbb{R}^N}E(u) = \sum_{i,j}\frac{\epsilon}{2}w_{ij}(u(i)-u(j))^2 +\frac{1}{\epsilon}\mathbb{W}(u) + \frac{\eta}{2}\sum_i\Lambda(i)(u(i)-y(i))^2
\end{equation}
where $u$ represents the labels of image or data and $i$ and $j$ are the index of data. The value of $\epsilon$ is typically larger than its counterpart in the phase field model, while $\eta$ is a parameter that balances the phase model and the prior term. $\mathbb{W}(\cdot)$ is obtained by summing a series of double-well functions, denoted as $\mathbb{W}(u) = \frac14\sum_{i=1}^{N}(u(i)^2 - 1)^2$. Here, the diagonal matrix $\Lambda$ and the vector $y$ are completely dependent on the prior and known labels \cite{Shensun2023}. 

Compared to the traditional phase-field model, the key distinction in the graphical model lies in the parameter $w_{ij}$, where $i$ and $j$ denote two data points. The weight is determined through two components: $w_{ij} =K(i,j) \cdot N(i,j)$, with $K(i,j)$ denoting a feature similarity factor and $N(i,j)$ representing a proximity indicator. The similarity factor is computed using a Gaussian kernel function: $K(i,j) = \exp\left(-\|P_{i} - P_{j}\|^2/\sigma^2\right)$, where $P_{i}$ stands for the feature of the data point $i$ and $\sigma^2$ is a parameter controlling the kernel width. The proximity term serves as a binary function that assesses data adjacency for long-distance interactions. We refer to see \cite[Section 2.2]{Shensun2023} for further information, including establishing a nonlocal and graph Laplacian.

For the image segmentation, we choose the model parameter as follows: $\epsilon = \eta = 10$. In the experiment described below, we examine using Algorithm \ref{alg:cap} to address a nonconvex and nonlinear problem. The energy function $E(u)$ is divided into two components: $H(u) =  \sum_{ij}\frac{\epsilon}{2}w_{ij}(u(i)-u(j))^2 + \frac{\eta}{2}\sum_i\Lambda(i)(u(i)-y(i))^2$ and $F(u) = \frac{1}{\epsilon}\mathbb{W}(u)$. Transitioning from $H(u) + F(u)$ to $H^n(u) - F^n(u)$ constitutes a special convex splitting approach at the $n$-th iteration. Moreover, the parameter $\delta t$ must satisfy $\delta t<\frac{2}{3L}$, where $L$ is the Lipschitz constant of $F(u)$. For other algorithms in Table \ref{tab:nonlocalGLmodel}, we use the following convex splitting:
\begin{align*}
    E(u) = \sum_{i,j}\frac{\epsilon}{2}w_{ij}(u(i)-u(j))^2 &+ \frac{\eta}{2}\sum_i\Lambda(i)(u(i)-y(i))^2 + \frac{L}{2}\sum_i u(i)^2 - \left(\frac{L}{2}\sum_i u(i)^2 +\frac{1}{\epsilon}\mathbb{W}(u)\right).
\end{align*}

Table \ref{tab:nonlocalGLmodel} provides a comparison of 9 algorithms: DCA \cite{LeThi2018}, pDCA$_e$ \cite[Section 3]{Wen2018} (or \cite{Liu2019}), BDCA \cite[Section 3]{Artacho2018}, as well as Algorithms \ref{alg:cap} including FlowBA (no preconditioners and line search), FlowBAp, FlowBApe, FlowBBA (no preconditioners), FlowBBAp, and FlowBBApe, which are variations of Algorithm \ref{alg:cap}. Specifically, DCA, FlowBA, pDCA$_e$, BDCA, and FlowBBA utilize the conjugate gradient (CG) method to solve the subproblem of the DC program with a precision requirement of $\|u^{l+1} - u^l\|<10^{-8}$, where $\{u^l\}_l$ is the generated sequence by solving the subproblem. On the other hand, algorithms FlowBAp, FlowBApe, FlowBBAp, and FlowBBApe employ a preconditioned scheme, which only uses 50 parallel preconditioned Jacobi iterations. The first criterion ``DICE Bound" in Table \ref{tab:nonlocalGLmodel}, relies on the segmentation coefficient DICE, a metric used to evaluate the segmentation outcome. It is calculated using the formula $\text{DICE} = \frac{2|X\cap Y|}{|X|+|Y|}$, where $|X|$ and $|Y|$ denote the pixel counts in the segmentation and ground truth, respectively. Moreover, the criterion ``DICE Bound" refers to the DICE coefficient of the segmentation, achieving a value of $0.993$, indicating a high degree of similarity.

Analyzing Table \ref{tab:nonlocalGLmodel}, it is evident that the algorithm FlowBBAp and FlowBBApe outperform others in terms of running time, requiring the least amount of time to fulfill termination criteria. Notably, the BDCA algorithm exhibits the fewest iterations among all algorithms, attributed to the high accuracy of its subproblem solution and the line search.
\begin{sidewaystable}[htbp!]
\caption{Solving \eqref{eq:nonlocalGL} on seven different termination criteria.(Criteria \uppercase\expandafter{\romannumeral1}: $\|\nabla E(u)\|$, Criteria \uppercase\expandafter{\romannumeral2}: $\|u^n-u^{n-1}\|$)\\
\protect\textcircled{1}DCA \protect\textcircled{2}FlowBA \protect\textcircled{3}FlowBAp \protect\textcircled{4}FlowBApe \protect\textcircled{5}pDCA$_{\text{e}}$ \protect\textcircled{6}BDCA \protect\textcircled{7}FlowBBA \protect\textcircled{8}FlowBBAp \protect\textcircled{9}FlowBBApe\label{tab:nonlocalGLmodel}}
\begin{tabular}{cccrrrrrrrrr}
\hline
\multicolumn{3}{c}{Criteria} &{\textcircled{1}}  &{\textcircled{2}} & {\textcircled{3}}  &{\textcircled{4}} & {\textcircled{5}}  &{\textcircled{6}} &{\textcircled{7}} & {\textcircled{8}}  &{\textcircled{9}}\\ \hline
\multicolumn{2}{c}{\multirow{2}{*}{DICE Bound}} & Iter & 16 & 27 & 36 & 33& 94 & \textbf{10} & 16 & 17 & 18\\ 
& &Time(s)& 57.90 & 60.00 & 21.93 & 22.38& 45.07 & 38.05 & 37.98 & \textbf{14.59} &15.31 \\ \hline
\multirow{6}{*}{\uppercase\expandafter{\romannumeral1}}& \multirow{2}{*}{$10^{-1}$}& Iter & 22 & 39 & 52 & 47& 129 & \textbf{14} & 27 & 32 & 32\\
& &Time(s) & 78.83 & 81.73 &31.19 &30.24 & 60.87 & 51.82 & 62.27 & \textbf{23.84} &23.96\\
 & \multirow{2}{*}{$10^{-3}$}&Iter & 40 & 52 & 69 & 63& 184 & \textbf{23} & 38 & 45 & 44\\
& &Time(s)& 142.84 & 104.64 & 41.04 & 39.66 & 86.32 & 84.24 & 87.03 & {31.87} & \textbf{31.38}\\ 
& \multirow{2}{*}{$10^{-5}$} & Iter & 77 & 98 & 130 & 119 & 301 & \textbf{42} & 55 & 66 &62 \\
& &Time(s)& 279.96 & 187.23 & 76.43 & 72.58 & 140.83 & 159.47 & 126.79 & {44.85} & \textbf{42.48} \\  \hline
\multirow{6}{*}{\uppercase\expandafter{\romannumeral2}} & \multirow{2}{*}{$10^{-1}$} &Iter& 24 & 41 & 54 & 49& 112 & \textbf{15} & 29 & 34 &33 \\
& &Time(s)& 85.89 & 85.23 & 32.34 &  31.42& 53.17 & 55.46 & 66.76 & {25.07} & \textbf{24.58}\\ 
& \multirow{2}{*}{$10^{-3}$}  &Iter& 50 & 62 & 78 & 72& 173 & \textbf{29} & 40 & 47 & 46 \\
& &Time(s)& 178.43 & 122.38 & 46.27 & 44.94 & 81.22 & 105.94 & 91.58 & {33.10} & \textbf{32.62}\\
& \multirow{2}{*}{$10^{-5}$}&Iter& 87 & 108 & 139 &128 & 268 & \textbf{47} & 60 & 71 & 72\\
& &Time(s)& 315.56 & 206.31 & 81.64 & 77.89 & 125.43 & 178.65 & 139.53 & \textbf{47.94} &49.45 \\ \hline
\end{tabular}
\end{sidewaystable}
\begin{figure}[htbp]   
    \caption{Segmentation results for the stone image}
    \label{fig:stone_seg}    
    \centering            
    \begin{subfigure}[b]{0.17\textwidth}
        \includegraphics[width=\linewidth]{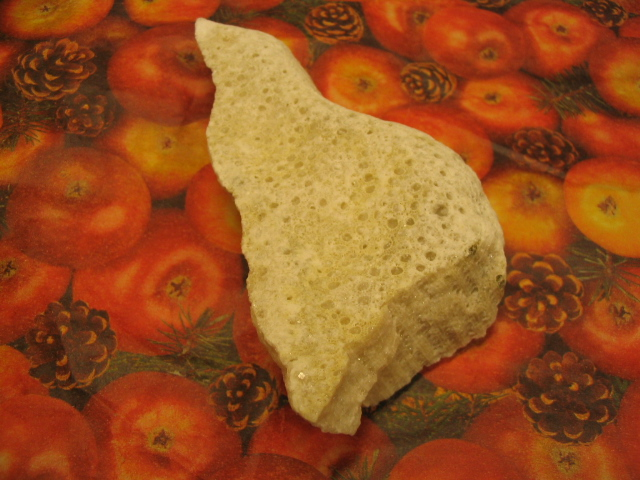}
        \caption{Image} 
        \label{fig:stone}
    \end{subfigure}\hfill
    \begin{subfigure}[b]{0.17\textwidth}
        \includegraphics[width=\linewidth]{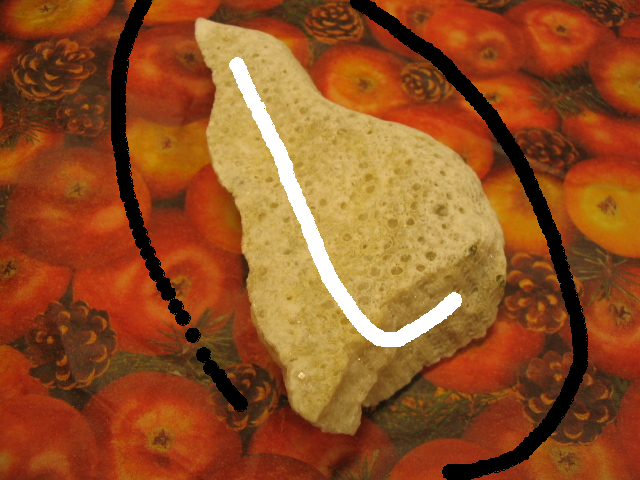}
        \caption{Label} 
        \label{fig:stonelable}
    \end{subfigure}\hfill
    \begin{subfigure}[b]{0.17\textwidth}
        \includegraphics[width=\linewidth]{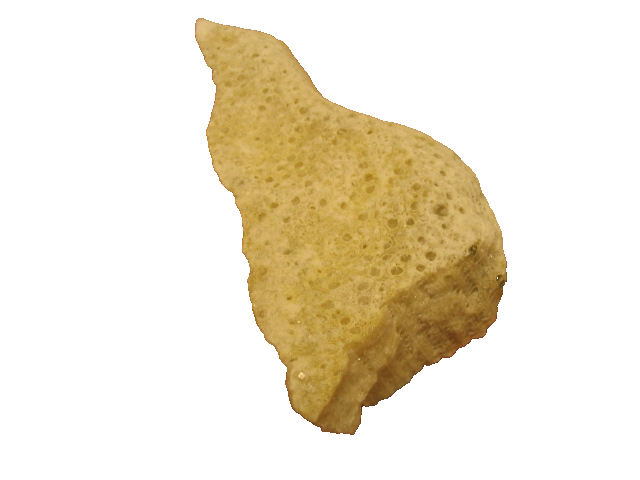}
        \caption{result1} 
        \label{fig:result1}
    \end{subfigure}\hfill
    \begin{subfigure}[b]{0.17\textwidth}
        \includegraphics[width=\linewidth]{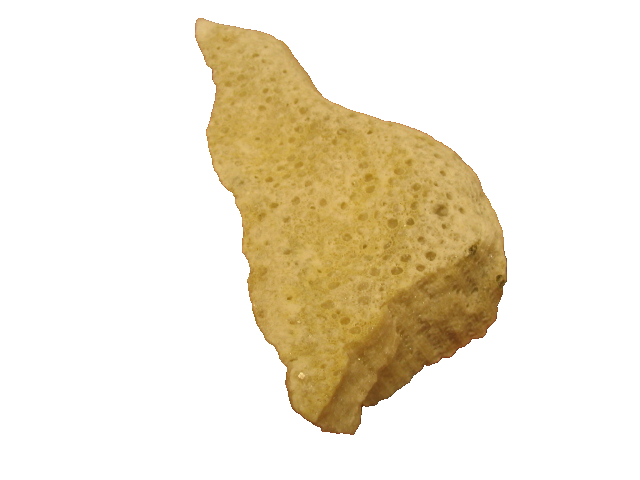}
        \caption{result2} 
        \label{fig:result2}
    \end{subfigure}\hfill
    \begin{subfigure}[b]{0.17\textwidth}
        \includegraphics[width=\linewidth]{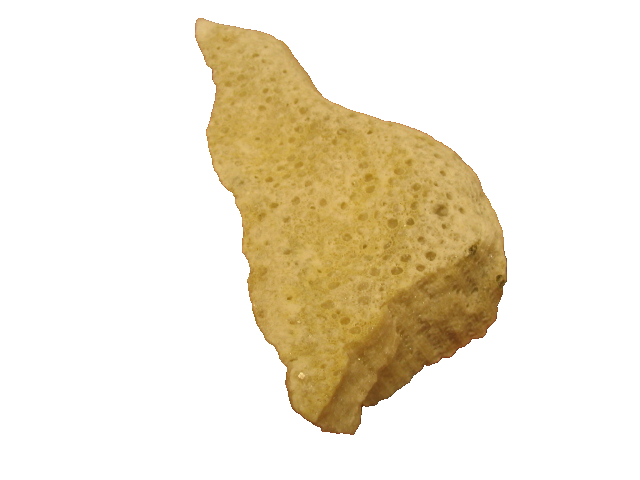}
        \caption{result3} 
        \label{fig:result3}
    \end{subfigure}
\end{figure}

In Figure \ref{fig:stone_seg}, we illustrate the segmentation result of an image. Figure \ref{fig:stone} is the original image, Figure \ref{fig:stonelable} labels the prior that the white pixels label the positive prior. The black pixels label the negative prior, and the last three figures represent the segmentation results under three different criteria (DICE Bound, $\|\nabla E(u)\|\leq 10^{-5}$ and $\|u^n-u^{n-1}\|\leq 10^{-5}$).

Let us end this section with a remark for the KL properties of the modified SCAD regularization \eqref{eq:lsp} and the graph Ginzburg-Landau functional \eqref{eq:nonlocalGL}.
\begin{remark}\label{rem:KL:conver:two:models}
The discrete graph Ginzburg-Landau functional, a polynomial function of $u$, is also a semi-algebraic function and possesses the KL property \cite[Section 2.2]{ABS}. Its corresponding semi-algebraic set can be written as
\begin{equation*}
    \left\{ (u, z) \;\middle|\; z - \sum_{ij} \frac{\epsilon}{2} w_{ij} (u(i) - u(j))^2 + \frac{1}{4\epsilon} \sum_i (u(i)^2 - 1)^2 + \frac{\eta}{2} \sum_i \Lambda(i) (u(i) - y(i))^2 = 0 \right\}.
\end{equation*}
Regarding the modified SCAD regularization \eqref{eq:huber-scad}, as $ p_{M}(u_i)$ represents a one-dimensional piecewise quadratic function, by similar arguments as in \cite[Section 5.2]{Li2018}, we conclude that $E(u)$ is a KL function. Furthermore, with Lemma \ref{lem:level_bound_A}, we confirm that $A(x,y)$ (or $\tilde A(x,y)$) is also a KL function and the global convergence of these two models thus follows.  
\end{remark}

\section{Conclusion}\label{sec:conclusion}
We proved the global convergence of the iterative sequence in finite-dimensional spaces for the commonly used second-order convex splitting method in phase-field simulations. By incorporating line search and preconditioning techniques, we improved its performance with global convergence. For the convergence of the continuous flow in Banach spaces, the \L ojasiewicz inequality in Banach spaces \cite{simon,haraux,bolte1} is a possible way. It would be very interesting to design efficient preconditioners for the Cahn-Hilliard equation \cite{Elliott1986, aasuli, Elliott1987, Elliott2021}. The efficacy of the proposed algorithm was validated through numerical experiments. Additionally, we find the third- to sixth-order IMEX splitting methods \cite{ARW1995} particularly intriguing within the context of the current line search and preconditioning framework. The lower-complexity bounds \cite{xuyangyang} are also relevant to the proposed algorithm. In addition, the non-monotone line search \cite{Dai2002, Ferreira2024} can be used when energy decay is not available. 

\bibliographystyle{plain}
\bibliography{bdfab_ls_JNA}

\end{document}